\documentclass[a4paper]{article}
\usepackage{RR}
\RRNo{6451}
\usepackage{hyperref}
\usepackage{color}

\RRdate{F\'evrier 2008}

\usepackage{amsmath}
\usepackage{amsfonts}
\graphicspath{{Figures/}{logos/}}
\newcommand\R{\mathbb{R}}
\newtheorem{remark}{Remark}
\newtheorem{proposition}{Proposition}
\RRauthor{
Jacques Sainte-Marie\thanks{Saint-Venant Laboratory, 6 quai Watier, 78400 Chatou}%
  \thanks[sfn]{INRIA Rocquencourt, B.P.~105, 78153 Le Chesnay Cedex, France.}%
  \and
Marie-Odile Bristeau
\thanksref{sfn}
}
\authorhead{Sainte-Marie \& Bristeau}
\RRtitle{D\'erivation d'un mod\`ele de type Saint-Venant non hydrostatique;\\ Comparaison
avec les mod\`eles de type Boussinesq}
\RRetitle{Derivation of a non-hydrostatic shallow water model;\\ Comparison with Saint-Venant and Boussinesq systems}
\titlehead{Derivation of a non-hydrostatic shallow water model}
\RRresume{A partir des \'equations de Navier-Stokes \`a surface libre,
on obtient deux mod\`eles moyenn\'es sur la verticale, non hydrostatiques
qui \'etendent le syst\`eme de Saint-Venant et incluent le frottement et
la viscosit\'e. La complexit\'e des
formulations obtenues d\'epend du niveau d'approximation retenu pour la
pression du fluide. Les mod\`eles obtenus sont compar\'es aux
formulations de type Boussinesq.
}
\RRabstract{From the free surface Navier-Stokes system, we derive the
 non-hydrostatic Saint-Venant
 system for the shallow waters including friction and viscosity. The derivation
leads to two formulations of growing complexity depending on the level
 of approximation chosen for the
 fluid pressure. The obtained models are compared with the Boussinesq models.}
\RRmotcle{Equations de Navier-Stokes, \'equations de Saint-Venant,
\'equations de Boussinesq, surface libre, termes dispersifs}
\RRkeyword{Navier-Stokes equations, Saint-Venant equations, Boussinesq equations, Free surface, Dispersive terms}
\RRprojets{MACS et BANG}
 \RRtheme{\THNum} 
\URRocq 

\begin{document}
\makeRR   

\tableofcontents

\section{Introduction}

Despite the available numerical results obtained by the simulation of
the Navier-Stokes equations, there exists a demand for models of reduced
complexity such as shallow waters type models.

\begin{figure}[htbp]
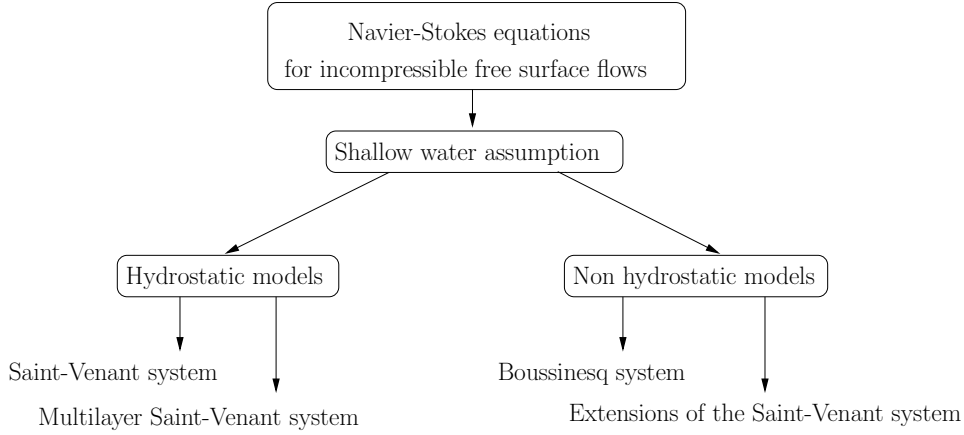

\begin{center}
\resizebox{13cm}{!}{\input Figures/models.pstex_t}
\caption{Averaged models derived from Navier-Stokes equations.}
\label{fig:shallow}
\end{center}
\end{figure}

Non-linear shallow water equations model the dynamics of a shallow,
rotating layer of homogeneous incompressible fluid and are typically
used to describe vertically averaged flows in two or three dimensional domains,
in terms of horizontal velocity and depth variation, see Fig.~\ref{fig:shallow}. This set of equations is particularly well-suited for the
study and numerical simulations of a large class of geophysical phenomena,
such as rivers, coastal domains, oceans, or even run-off or avalanches when
modified with adapted source terms \cite{bouchut}.

The classical Saint-Venant system \cite{saint-venant}
with viscosity and friction \cite{gerbeau,marche,saleri} is well suited
for modeling of dam breaks or hydraulic jump but due to the hydrostatic
assumption it is not well adapted for the modeling of gravity waves propagation.

For the modeling of long wavelength, small amplitude, gravity waves, the
Boussinesq system \cite{boussinesq1,boussinesq2,boussinesq3} is
used. The Boussinesq equations are obtained from the Euler equations
i.e. ignoring rotational and dissipative effects
\cite{bbm,cienfuegos1,cienfuegos2,nwogu,peregrine,walkley}. In practice, the use of such models
ignoring rotational and friction effects at the bottom may be very
restrictive. Furthermore, even when well posed, the Boussinesq models
often exhibit a lack of conservation energy that is odd since they are
derived from Euler equations \cite{saut1,saut2}.

The objective of this paper is twofold. First, we want to extend the
Saint-Venant system so that the long waves propagation can be modeled
and second we aim at comparing/unifying the obtained formulation with
the Boussinesq system, see Fig.~\ref{fig:shallow}. The paper is organized as follows. In section
\ref{sec:NS}, we recall the Navier-Stokes system with a free moving
boundary and its closure. We also present the Saint-Venant and
Boussinesq assumptions and the associated rescaling. In section
\ref{sec:SW} we recall the Shallow Water system and show the hydrostatic
Boussinesq system assumption corresponds to the classical Saint-Venant
system. In section \ref{sec:SV_nhyd1}, the
hydrostatic assumption is relaxed and we obtain two formulations of
growing complexity extending the Saint-Venant system and depending on the level of approximation chosen for the fluid pressure.

\section{The Navier-Stokes system}
\label{sec:NS}

\begin{figure}[htbp]
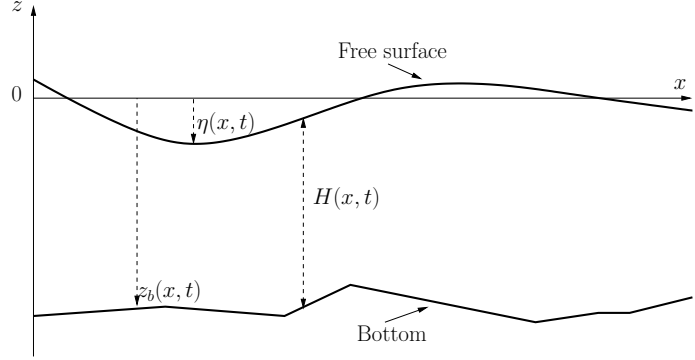

\begin{center}
\resizebox{9cm}{!}{\input Figures/notations.pstex_t}
\caption{Notations: water height $H(x,t)$, free surface $\eta(x,t)$ and bottom $z_b(x,t)$.}
\label{fig:notations}
\end{center}
\end{figure}

Let start with the Navier-Stokes system \cite{lions} restricted to two dimensions
with gravity in which the $z$ axis represents the vertical
direction. For simplicity, the viscosity will be kept constant
throughout the paper. Therefore we have the following general formulation
expression:
\begin{eqnarray}
& & \frac{\partial u}{\partial x} + \frac{\partial w}{\partial z} = 0,\label{eq:NS_2d1}\\
& & \frac{\partial u}{\partial t} + u \frac{\partial u}{\partial x} + w \frac{\partial u}{\partial z} + \frac{1}{\rho} \frac{\partial p}{\partial x} = \frac{\partial \Sigma_{xx}}{\partial x} + \frac{\partial \Sigma_{xz}}{\partial z},\label{eq:NS_2d2}\\
& & \frac{\partial w}{\partial t} + u\frac{\partial w}{\partial x} + w\frac{\partial w}{\partial z} + \frac{1}{\rho} \frac{\partial p}{\partial z} = -g + \frac{\partial \Sigma_{zx}}{\partial x} + \frac{\partial \Sigma_{zz}}{\partial z},
\label{eq:NS_2d3}
\end{eqnarray}
and we consider this system for
$$t>t_0, \quad x \in \R, \quad z_b(x,t) \leq z \leq \eta(x,t),$$
where $\eta(x,t)$ represents the free surface elevation, ${\bf u}=(u,w)^T$ the horizontal and vertical velocities. The water
height is $H = \eta - z_b$, see Fig.~\ref{fig:notations}. We consider the
bathymetry $z_b$ can vary with respect to abscissa $x$ and also with
respect to time $t$. The chosen form of the viscosity tensor is
\begin{eqnarray*}
 \Sigma_{xx} = 2 \nu \frac{\partial u}{\partial x}, & & \Sigma_{xz} = \nu \bigl( \frac{\partial u}{\partial z} + \frac{\partial w}{\partial x} \bigr),\label{eq:visco1}\\
\Sigma_{zz} = 2 \nu \frac{\partial w}{\partial z}, && \Sigma_{zx} = \nu \bigl( \frac{\partial u}{\partial z} + \frac{\partial w}{\partial x}\bigr),
\label{eq:visco2}
\end{eqnarray*}
with $\nu$ the viscosity coefficient. For a more complex form of the
viscosity tensor using eddy and bulk viscosities, the reader can refer
to \cite{levermore}.

\subsection{Boundary conditions}

The system (\ref{eq:NS_2d1})-(\ref{eq:NS_2d3}) is complete with boundary conditions. The
outward and upward unit normals to the free surface ${\bf n}_s$ and to the bottom ${\bf n}_b$ are
given by
$${\bf n}_s = \frac{1}{\sqrt{1 + \bigl(\frac{\partial \eta}{\partial x}\bigr)^2}}  \left(\begin{array}{c} -\frac{\partial \eta}{\partial x}\\ 1 \end{array} \right), \quad {\bf n}_b = \frac{1}{\sqrt{1 + \bigl(\frac{\partial z_b}{\partial x}\bigr)^2}}  \left(\begin{array}{c} -\frac{\partial z_b}{\partial x}\\ 1 \end{array} \right).$$
Let
$\Sigma_T$ be the total stress tensor with
$$\Sigma_T = -\frac{1}{\rho} p I_d + \left(\begin{array}{cc} \Sigma_{xx} & \Sigma_{xz}\\ \Sigma_{zx} & \Sigma_{zz}\end{array}\right).$$

\subsubsection{At the free surface}

Classically at the free surface we have the kinematic boundary condition
\begin{equation}
\frac{\partial \eta}{\partial t} + u_s \frac{\partial \eta}{\partial x}
-w_s = 0,
\label{eq:free_surf} 
\end{equation}
where the subscript $s$ denotes the value of the
considered quantity at the free surface. Considering the air viscosity is
negligible, the continuity
of stresses at the free boundary imposes
\begin{equation}
\Sigma_T {\bf n}_s = -\frac{p^a}{\rho} {\bf n}_s,
\label{eq:BC_h}
\end{equation}
where $p^a=p^a(x,t)$ is a given function corresponding to the
atmospheric pressure. Relation~(\ref{eq:BC_h}) is
equivalent to
$${\bf n}_s . \Sigma_T {\bf n}_s = -\frac{p^a}{\rho},\quad \hbox{and}\quad {\bf t}_s  .\Sigma_T {\bf n}_s = 0,$$
${\bf t}_s$ being orthogonal to ${\bf n}_s$.

\subsubsection{At the bottom}

Since we consider the bottom can vary with respect to time $t$, the
kinematic boundary condition is
\begin{equation}
\frac{\partial z_b}{\partial t} + u_b \frac{\partial z_b}{\partial x}
-w_b = 0,
\label{eq:bottom} 
\end{equation}
where the subscript $b$ denotes the value of the
considered quantity at the bottom and $(x,t) \mapsto z_b(x,t)$ is a
given function. Note that Eq.~(\ref{eq:bottom})
reduces to a classical no-penetration condition when $z_b$
does not depend on time $t$.

For the stresses at the bottom we consider a wall law under the form
\begin{equation}
\Sigma_T {\bf n}_b - ({\bf n}_b . \Sigma_T {\bf n}_b){\bf n}_b = \kappa({\bf v_b},H) {\bf v}_b,
\label{eq:BC_z_b}
\end{equation}
with ${\bf v}_b={\bf u}_b - (0,\frac{\partial z_b}{\partial t})^T$ the
relative velocity between the water and the bottom. If $\kappa({\bf v_b},H)$ is constant then
we recover a Navier friction condition as in \cite{gerbeau}. Introducing
laminar $k_l$ and turbulent $k_t$ friction, we use the expression 
$$\kappa({\bf v_b},H,\nu) = k_l + k_t H |{\bf v_b}|,$$
corresponding to the boundary condition used in \cite{marche}. Another
form of $\kappa({\bf v_b},H)$ is used in \cite{bouchut} and for other
wall laws, the reader can also refer
to \cite{valentin}. Due to thermomechanical considerations, in the
sequel we suppose $\kappa({\bf v_b},H) \geq 0$ and $\kappa({\bf
v_b},H)$ is often simply denoted $\kappa$.

Let ${\bf t}_b$ satisfying ${\bf t}_b . {\bf n}_b = 0$
then when multiplied by ${\bf t}_b$ and ${\bf n}_b$, Eq.~(\ref{eq:BC_z_b}) leads to
$${\bf t}_b. \Sigma_T {\bf n}_b = \kappa {\bf v}_b . {\bf t}_b,\quad \hbox{and}\quad {\bf v}_b. {\bf n}_b = 0.$$
\begin{remark}
If the boundary condition (\ref{eq:BC_z_b}) was written under the form $\Sigma_T. n_b = \kappa({\bf v_b},H) {\bf v}_b$ as in Ferrari {\it et al.}~\cite[Eq.~(2.25), p.~217]{saleri}, then in absence of friction and viscosity, this would give
$$p_b=0,$$
that is not correct.
\label{rem:BC_gravity}
\end{remark}

\subsection{The rescaled system}
\label{subsec:NS_rescaled}

The physical system is rescaled using the quantities
\begin{itemize}
\item $h$ and $\lambda$, two characteristic dimensions along the $z$ and $x$
      axis respectively,
\item $a_s$ the typical wave amplitude, $a_b$ the typical bathymetry variation,
\item $C=\sqrt{gh}$ the typical horizontal wave speed.
\end{itemize}
Classically for the derivation of the Saint-Venant system, we introduce
the small parameter
$$\varepsilon = \frac{h}{\lambda}.$$
When considering long waves propagation,
another important parameter needs be considered, namely
$$\delta = \frac{a_s}{h},$$
and we consider for the bathymetry $\frac{a_b}{h} = {\cal O}(\delta)$.
Depending on the application, $\delta$ can be considered or not as a small
parameter. For finite amplitude wave theory and assuming
$z_b(x,t)=z_b^0$, one considers $\varepsilon \ll 1$, $\delta =
{\cal O}(1)$ whereas the Boussinesq waves theory requires
$$\delta \ll 1,\quad\varepsilon \ll 1\quad\hbox{and}\quad U_r = {\cal O}(1).$$
where $U_r$ is the Ursell number defined by $U_r = \frac{\delta}{\varepsilon^2}$, see \cite{ursell}. All along this work, we consider $\varepsilon \ll 1$ whereas, even if the parameter $\delta$ is introduced in the
rescaling, the assumption $\delta \ll 1$ is not considered
(paragraphs~\ref{subsec:extension_1}, \ref{subsec:nhyd1} and \ref{subsec:prop}) except when
explictly mentioned.

As for the Saint-Venant system
\cite{gerbeau,marche}, we introduce some characteristic quantities~:
$T=\lambda/C$ for the time, $W = a_s/T = a_b/T = \varepsilon\delta C$
for the vertical velocity, $U = W/\varepsilon=\delta C$, for the horizontal velocity, $P=\rho C^2$ for
the pressure. This leads to the following dimensionless quantities
$$\tilde{x} = \frac{x}{\lambda},\quad \tilde{z} = \frac{z}{h},\quad\tilde{\eta} = \frac{\eta}{a_s},\quad\tilde{t} = \frac{t}{T},$$
$$\tilde{p} = \frac{p}{P},\quad \tilde{u} = \frac{u}{
U},\quad\mbox{and}\quad \tilde{w} = \frac{w}{W}.$$
Note that the definition of the charateristic velocites implies $\delta = \frac{U}{C}$ so $\delta$ also corresponds to the Froude number. When $\delta = {\cal O}(1)$ we have $U \approx C$ and we recover the classical rescaling used for the Saint-Venant system. For the bathymetry $z_b$ we write $z_b(x,t) = Z_b(x)+b(t)$ and we
introduce $\tilde{z}_b = Z_b/h$ and $\tilde{b} = b / a_b$. This leads to
$$\frac{\partial z_b}{\partial t} = \varepsilon\delta C
\frac{\partial\tilde{b}}{\partial \tilde{t}} = W\frac{\partial\tilde{b}}{\partial \tilde{t}},\quad\mbox{and}\quad\frac{\partial z_b}{\partial x} = \varepsilon\frac{\partial\tilde{z}_b}{\partial \tilde{x}}.$$
The different rescaling applied to the time
and space derivatives of $z_b$ means that a classical shallow water
assumption is made concerning the space variations of the bottom profile whereas we assume
the time variations of $z_b$ lie in the framework of the Boussinesq
assumption and are consistent with the rescaling applied to the velocity
$w$.

We also introduce $\tilde{\nu} = \frac{\nu}{\lambda C}$ and we set $\tilde{\kappa} = \frac{\kappa}{C}$.
Note that the definitions for the dimensionless quantities are
consistent with the one used for the Boussinesq
system \cite{peregrine,walkley}. Note also that the rescaling used by
Nwogu \cite{nwogu} differs from the preceding one since Nwogu uses
$\tilde{w} = \frac{\varepsilon^2}{W}w$.

As in \cite{gerbeau,marche}, we suppose we are in the following asymptotic regime
$$\tilde{\nu} = \varepsilon \nu_0,\qquad \mbox{and}\qquad\tilde{\kappa} =
\varepsilon \kappa_0,$$
with $\kappa_0 = \kappa_{l,0} + \varepsilon\kappa_{t,0}(\tilde{\bf
v}_b,\tilde{H})$, $\kappa_{l,0}$ being constant.

This non-dimensionalization of the system (\ref{eq:NS_2d1})-(\ref{eq:NS_2d3}) leads to
\begin{eqnarray}
\lefteqn{\frac{\partial \tilde{u}}{\partial \tilde{x}} + \frac{\partial \tilde{w}}{\partial \tilde{z}} = 0,\label{eq:NS_2d_re1}}\\
\lefteqn{\varepsilon\delta\frac{\partial \tilde{u}}{\partial \tilde{t}} + \varepsilon\delta^2\tilde{u} \frac{\partial \tilde{u}}{\partial \tilde{x}} + \varepsilon\delta^2\tilde{w} \frac{\partial \tilde{u}}{\partial \tilde{z}} + \varepsilon\frac{\partial \tilde{p}}{\partial \tilde{x}} = \varepsilon^2\delta\frac{\partial }{\partial \tilde{x}}\left(2\nu_0\frac{\partial \tilde{u}}{\partial \tilde{x}}\right) \nonumber}\\
& & + \frac{\partial }{\partial \tilde{z}}\left( \delta\nu_0\frac{\partial \tilde{u}}{\partial \tilde{z}} + \varepsilon^2\delta\nu_0\frac{\partial \tilde{w}}{\partial \tilde{x}}\right), \label{eq:NS_2d_re2}\\
\lefteqn{\varepsilon^2\delta\left(\frac{\partial \tilde{w}}{\partial \tilde{t}} + \delta\tilde{u}\frac{\partial \tilde{w}}{\partial \tilde{x}} + \delta\tilde{w}\frac{\partial \tilde{w}}{\partial \tilde{z}} \right) + \frac{\partial \tilde{p}}{\partial \tilde{z}} = -1 \nonumber}\\
& & + \frac{\partial }{\partial \tilde{x}}\left( \varepsilon\delta\nu_0\frac{\partial \tilde{u}}{\partial \tilde{z}} + \nu_0\varepsilon^3\delta \frac{\partial \tilde{w}}{\partial \tilde{x}}\right)
 + \varepsilon\delta\frac{\partial }{\partial \tilde{z}}\left(2\nu_0\frac{\partial \tilde{w}}{\partial \tilde{z}}\right), \label{eq:NS_2d_re3}
\end{eqnarray}
with the boundary conditions (\ref{eq:free_surf}), (\ref{eq:BC_h}), (\ref{eq:bottom})  and (\ref{eq:BC_z_b}) becoming
\begin{eqnarray}
\lefteqn{\frac{\partial \tilde{\eta}}{\partial \tilde{t}} + \delta\tilde{u}_s \frac{\partial \tilde{\eta}}{\partial \tilde{x}} - \tilde{w}_s = 0,\label{eq:BC_re1}}\\
\lefteqn{2\varepsilon\delta\nu_0\left.\frac{\partial \tilde{w}}{\partial \tilde{z}}\right|_s -\tilde{p}_s -\varepsilon^2\delta^2\nu_0\frac{\partial \tilde{\eta}}{\partial \tilde{x}}\left(\left.\frac{\partial \tilde{u}}{\partial \tilde{z}}\right|_s + \varepsilon^2\left.\frac{\partial \tilde{w}}{\partial \tilde{x}}\right|_s\right) = -\delta\tilde{p}^a,\label{eq:BC_re2}}\\
\lefteqn{\delta\nu_0\left(\left.\frac{\partial \tilde{u}}{\partial \tilde{z}}\right|_s + \varepsilon^2\left.\frac{\partial \tilde{w}}{\partial \tilde{x}}\right|_s \right) - \varepsilon\delta\frac{\partial \tilde{\eta}}{\partial \tilde{x}} \left(2\varepsilon\delta\nu_0\left.\frac{\partial \tilde{u}}{\partial \tilde{x}}\right|_s - \tilde{p}_s\right) = \varepsilon\delta^2\frac{\partial \tilde{\eta}}{\partial \tilde{x}} \tilde{p}^a, \label{eq:BC_re3}}\\
\lefteqn{\frac{\partial \tilde{b}}{\partial \tilde{t}} + \tilde{u}_b \frac{\partial \tilde{z}_b}{\partial \tilde{x}} - \tilde{w}_b = 0, \label{eq:BC_re4}}\\
\lefteqn{\delta\nu_0\left(\varepsilon^2\left.\frac{\partial \tilde{w}}{\partial \tilde{x}}\right|_b + \left.\frac{\partial \tilde{u}}{\partial \tilde{z}}\right|_b \right) - \varepsilon\frac{\partial \tilde{z}_b}{\partial \tilde{x}} \left(2\varepsilon\delta\nu_0 \left.\frac{\partial \tilde{u}}{\partial \tilde{x}}\right|_b - p_b \right) \nonumber}\\
& & \quad + \varepsilon\frac{\partial \tilde{z}_b}{\partial \tilde{x}}\left(2\varepsilon\delta\nu_0\left.\frac{\partial \tilde{w}}{\partial \tilde{z}}\right|_b\right. - p_b - \left.\varepsilon\nu_0\frac{\partial \tilde{z}_b}{\partial \tilde{x}} \left(\delta\left.\frac{\partial \tilde{u}}{\partial \tilde{z}}\right|_b + \varepsilon^2\delta\left.\frac{\partial \tilde{w}}{\partial \tilde{x}}\right|_b\right)\right)\nonumber\\
& & = \varepsilon\delta\kappa_0 \sqrt{1 + \varepsilon^2\left(\frac{\partial \tilde{z}_b}{\partial \tilde{x}}\right)^2} \left(\tilde{u}_b + \varepsilon^2\frac{\partial \tilde{z}_b}{\partial \tilde{x}}\bigl(\tilde{w}_b - \frac{\partial \tilde{f}}{\partial \tilde{t}}\bigr)\right) \label{eq:BC_re5}.
\end{eqnarray}
For the sake of clarity, in the sequel we drop the symbol $\tilde{}$ and
we denote $\frac{\partial b}{\partial t} = \frac{\partial z_b}{\partial t}$.

\section{The Shallow Water system}
\label{sec:SW}

In this section we first derive the expression of the fluid pressure $p$ in
the context of the Shallow Water assumption and then show the
combination of the Boussinesq and hydrostatic assumption leads to the
classical Saint-Venant system.

The process used hereafter is similar to the technique employed by
Gerbeau and Perthame \cite{gerbeau} to derive a formulation for the viscous
Saint-Venant system.

\subsection{The vertically averaged system}

Using the divergence free condition, the system (\ref{eq:NS_2d_re1})-(\ref{eq:NS_2d_re3}) is
rewritten under the form
\begin{eqnarray}
\lefteqn{\frac{\partial u}{\partial x} + \frac{\partial w}{\partial z} = 0, \label{eq:NS_2d_newa1}}\\
\lefteqn{\varepsilon\delta\frac{\partial u}{\partial t} + \varepsilon\delta^2 \frac{\partial u^2}{\partial x} + \varepsilon\delta^2 \frac{\partial u w}{\partial z} + \varepsilon\frac{\partial p}{\partial x} = \varepsilon^2\delta\frac{\partial }{\partial x}\left(2\nu_0\frac{\partial u}{\partial x}\right)\nonumber}\\
& & \hspace*{1.5cm} + \frac{\partial }{\partial z}\left( \delta\nu_0\frac{\partial u}{\partial z} + \varepsilon^2\delta\nu_0 \frac{\partial w}{\partial x}\right),\label{eq:NS_2d_newa2}\\
\lefteqn{\varepsilon^2\delta\left(\frac{\partial w}{\partial t} + \delta \frac{\partial u w}{\partial x} + \delta\frac{\partial w^2}{\partial z} \right) + \frac{\partial p}{\partial z} = -1 \nonumber}\\
& & \hspace*{1.5cm} + \frac{\partial }{\partial x}\left( \varepsilon\delta\nu_0\frac{\partial u}{\partial z} + \varepsilon^3\delta\nu_0 \frac{\partial w}{\partial x}\right) + \varepsilon\delta\frac{\partial }{\partial z}\left(2\nu_0\frac{\partial w}{\partial z}\right).\label{eq:NS_2d_newa3}
\end{eqnarray}
Due to the applied rescaling some terms of the viscosity tensor e.g.
$$\varepsilon^3\delta\frac{\partial }{\partial x}\left( \nu_0 \frac{\partial w}{\partial x}\right)$$
are very small and could be neglected. But, as mentioned in
\cite[Remarks~1 and 2]{audusse}, the approximation of the viscous terms
have to preserve the dissipation energy that is an essential property of
the Navier-Stokes and averaged Navier-Stokes equations. Since we
privilege this stability requirement and in order to keep a symmetric
form of the viscosity tensor, we consider in the sequel a modified version of (\ref{eq:NS_2d_newa1})-(\ref{eq:NS_2d_newa3}) under the form
\begin{eqnarray}
\lefteqn{\frac{\partial u}{\partial x} + \frac{\partial w}{\partial z} = 0, \label{eq:NS_2d_new1}}\\
\lefteqn{\varepsilon\delta\frac{\partial u}{\partial t} + \varepsilon\delta^2 \frac{\partial u^2}{\partial x} + \varepsilon\delta^2 \frac{\partial u w}{\partial z} + \varepsilon\frac{\partial p}{\partial x} = \varepsilon^2\delta\frac{\partial }{\partial x}\left(2\nu_0\frac{\partial u}{\partial x}\right)\nonumber}\\
& & \hspace*{6cm} + \frac{\partial }{\partial z}\left( \delta\nu_0\frac{\partial u}{\partial z}\right),\label{eq:NS_2d_new2}\\
\lefteqn{\varepsilon^2\delta\left(\frac{\partial w}{\partial t} + \delta \frac{\partial u w}{\partial x} + \delta\frac{\partial w^2}{\partial z} \right) + \frac{\partial p}{\partial z} = -1 + \frac{\partial }{\partial x}\left( \varepsilon\delta\nu_0\frac{\partial u}{\partial z}\right)\nonumber}\\
& & \hspace*{6cm} + \frac{\partial }{\partial z}\left( 2\varepsilon\delta\nu_0\frac{\partial w}{\partial z}\right),\label{eq:NS_2d_new3}
\end{eqnarray}
corresponding to a viscosity tensor of the form
$$\Sigma_{xx} = 2 \nu \frac{\partial u}{\partial x}, \quad \Sigma_{xz} =
\Sigma_{zx} = \nu\frac{\partial u}{\partial z}, \quad \Sigma_{zz} = 2
\nu \frac{\partial w}{\partial z}.$$
\begin{remark}
If we strictly follow Audusse \cite[Lemma~2.1]{audusse}, the chosen form
 of the viscosity tensor will not allow us to include under the form of
 a square term in the energy equality the quantity
$$\nu_0\frac{\partial u}{\partial z}\frac{\partial w}{\partial x}.$$
But we will see in paragraph \ref{subsec:prop} that due to the shallow
 water assumption, this quantity appear as a friction term.
\end{remark}
From Eqs.~(\ref{eq:BC_re2}), it comes
$$p_s = \delta p^a + 2\varepsilon\delta\left.\frac{\partial w}{\partial z}\right|_s + {\cal O}(\varepsilon^2\delta^2),$$
so using Eqs.~(\ref{eq:BC_re3}) and (\ref{eq:BC_re5}) one obtains
\begin{equation}
\left.\frac{\partial u}{\partial z}\right|_s = {\cal O}(\varepsilon^2),\quad \left.\frac{\partial u}{\partial z}\right|_b = {\cal O}(\varepsilon),
\label{eq:BC1_ubis}
\end{equation}
and an integration of Eq.~(\ref{eq:NS_2d_new3}) from
$\delta\eta$ to $z$ gives
\begin{equation}
p - \delta p^a = \delta\eta - z + {\cal O}(\varepsilon\delta),
\label{eq:p_simple}
\end{equation}
leading to
$$\frac{\partial p}{\partial x} = {\cal O}(\delta).$$
The preceding relation inserted in (\ref{eq:NS_2d_new2}) leads to
\begin{equation}
\nu_0\frac{\partial^2 u}{\partial z^2} = {\cal O}(\varepsilon),
\label{eq:BC1_uter}
\end{equation}
and Eqs.~(\ref{eq:BC1_ubis}) and (\ref{eq:BC1_uter}) mean that
$$u(x,z,t) = u_0(x,t) + {\cal O}(\varepsilon),$$
i.e. we recognize the so-called ``motion by slices'' of the usual
Saint-Venant system. Then we introduce the averaged quantities
$$\bar{u} = \frac{1}{\delta \eta - z_b} \int_{z_b}^{\delta \eta} u\
dz,\quad \overline{u^2} = \frac{1}{\delta \eta - z_b} \int_{z_b}^{\delta \eta} u^2\ dz,$$
and the previous definitions involve
\begin{equation}
u(x,z,t) = \bar{u} + {\cal O}(\varepsilon),\quad\mbox{and}\quad \overline{u^2} = \bar{u}^2(x,z,t) + {\cal O}(\varepsilon).
\label{eq:u2}
\end{equation}
Note that the velocity $\bar{u}$ is exactly the one arising in the
conservation law for the water height since an integration of
Eq.~(\ref{eq:NS_2d_new1}) from $z_b$ to $\delta\eta$ with
boundary conditions (\ref{eq:BC_re1}) and (\ref{eq:BC_re4}) leads to
\begin{equation}
\frac{\partial \eta}{\partial t} - \frac{\partial z_b}{\partial t} + \frac{\partial}{\partial x} \bigl(H_\delta\bar{u}\bigr) = 0,
\label{eq:conserv}
\end{equation}
with $H_\delta = \delta\eta - z_b$. Conversely an integration of
Eq.~(\ref{eq:NS_2d_new1}) from $z_b$ to $z$ with
boundary conditions (\ref{eq:BC_re1}) and (\ref{eq:BC_re4}) leads to
\begin{equation}
w = \frac{\partial z_b}{\partial t} - \frac{\partial}{\partial
x}\int_{z_b}^z u\ dz = \frac{\partial z_b}{\partial t} - z\frac{\partial \bar{u}}{\partial x} +
\frac{\partial (z_b \bar{u})}{\partial x} + {\cal O}(\varepsilon).
\label{eq:cons_w}
\end{equation}
We use the approximations obtained in this paragraph to simplify the
boundary conditions (\ref{eq:BC_re1})-(\ref{eq:BC_re5}) and retaining
only the high order terms we obtain
\begin{eqnarray}
& &  \frac{\partial \eta}{\partial t} + \delta u_s \frac{\partial \eta}{\partial x} - w_s = 0,\label{eq:BC_ap1}\\
& & p_s = \delta p^a + 2\varepsilon\delta\nu_0\left.\frac{\partial w}{\partial z}\right|_s + {\cal O}(\varepsilon^3\delta),\label{eq:BC_ap2}\\
& & \delta\nu_0\left.\frac{\partial u}{\partial z}\right|_s - \varepsilon\delta\frac{\partial \eta}{\partial x} \left( 2\varepsilon\delta\nu_0\left.\frac{\partial u}{\partial x}\right|_s - p_s \right) = \varepsilon\delta^2\frac{\partial \eta}{\partial x} p^a, \label{eq:BC_ap3}\\
& & \frac{\partial z_b}{\partial t} + u_b \frac{\partial z_b}{\partial x} - w_b = 0, \label{eq:BC_ap4}\\
& & \delta\nu_0\left.\frac{\partial u}{\partial z}\right|_b - \varepsilon \frac{\partial z_b}{\partial x} \left( 2\varepsilon\delta\nu_0 \left.\frac{\partial u}{\partial x}\right|_b - p_b \right) = -\varepsilon\frac{\partial z_b}{\partial x} \left( 2\varepsilon\delta\nu_0\left.\frac{\partial w}{\partial z}\right|_b - p_b \right. \nonumber\\
& & \qquad \left. - \varepsilon\delta\nu_0\frac{\partial z_b}{\partial x}\left.\frac{\partial u}{\partial z}\right|_b \right) + \varepsilon\delta\kappa_0 \left(1 +
\frac{3\varepsilon^2}{2}\left(\frac{\partial z_b}{\partial x}\right)^2 \right) u_b + {\cal O}(\varepsilon^4\delta).\qquad\label{eq:BC_ap5}
\end{eqnarray}
Using the Leibniz rule i.e.
$$\frac{\partial}{\partial x} \int_{a(x)}^{b(x)} g\ dx_1 =
\int_{a(x)}^{b(x)} \frac{\partial g}{\partial x_1} dx_1 +
\frac{\partial b}{\partial x}g(a(x)) - \frac{\partial a}{\partial x}g(b(x)),$$
and the kinematic boundary conditions (\ref{eq:BC_ap1}) and (\ref{eq:BC_ap4}), an
integration of Eq.~(\ref{eq:NS_2d_new2}) from $z_b$ to $\delta\eta$
shows that a solution to (\ref{eq:NS_2d_new1})-(\ref{eq:NS_2d_new3}) satisfies
\begin{eqnarray*}
\varepsilon\delta\frac{\partial}{\partial t}\int_{z_b}^{\delta \eta} u\ dz + \varepsilon\frac{\partial}{\partial x}\int_{z_b}^{\delta \eta} \left(\delta^2 u^2 + p\right) dz = \varepsilon^2\delta\frac{\partial }{\partial x}\int_{z_b}^{\delta \eta} 2\nu_0\frac{\partial u}{\partial x} dz\nonumber\\
+ \delta\nu_0\left.\frac{\partial u}{\partial z}\right|_s - \varepsilon\delta \frac{\partial \eta}{\partial x} \left(2\varepsilon\delta\nu_0\left.\frac{\partial u}{\partial x}\right|_s - p_s \right)\nonumber\\
- \delta\nu_0\left.\frac{\partial u}{\partial z}\right|_b + \varepsilon \frac{\partial z_b}{\partial x} \left(2\varepsilon\delta\nu_0\left.\frac{\partial u}{\partial x}\right|_b - p_b \right),
\end{eqnarray*}
and using Eqs.~(\ref{eq:BC_ap3}) and (\ref{eq:BC_ap5}), we obtain
\begin{eqnarray}
\delta\frac{\partial}{\partial t}\int_{z_b}^{\delta \eta} u\ dz + \frac{\partial}{\partial x}\int_{z_b}^{\delta \eta} \left(\delta^2 u^2 + p\right) dz =  \varepsilon\delta\frac{\partial }{\partial x}\int_{z_b}^{\delta \eta} 2\nu_0\frac{\partial u}{\partial x} dz\nonumber\\
+ \delta^2 \frac{\partial \eta}{\partial x} p^a + \frac{\partial z_b}{\partial x} \left( 2\varepsilon\delta\nu_0\left.\frac{\partial w}{\partial z}\right|_b - p_b - \varepsilon\delta\nu_0\frac{\partial z_b}{\partial x}\left.\frac{\partial u}{\partial z}\right|_b \right) \nonumber\\
- \delta\kappa_0 \left( 1 + \frac{3\varepsilon^2}{2}\left(\frac{\partial z_b}{\partial x}\right)^2\right) u_b + {\cal O}(\varepsilon^3\delta),
\label{eq:momentum}
\end{eqnarray}
An expression for the pressure $p$ can be obtained as follows. An
integration of Eq.~(\ref{eq:NS_2d_new3}) from $z$ to $\delta \eta$ gives
\begin{eqnarray*}
\varepsilon^2\delta \int_z^{\delta \eta} \bigl(\frac{\partial w}{\partial t} + \delta\frac{\partial (uw)}{\partial x}\bigr) dz + \varepsilon^2\delta^2(w_s^2 - w^2) +
p_s - p = -(\delta \eta - z)\\
+ \varepsilon\delta\int_z^{\delta \eta} \frac{\partial}{\partial x}\left(\nu_0 \frac{\partial u}{\partial
z}\right)dz - 2\varepsilon\delta\nu_0\frac{\partial w}{\partial z} + \left.2\varepsilon\delta\nu_0\frac{\partial w}{\partial z}\right|_s
\end{eqnarray*}
and using the boundary conditions (\ref{eq:BC_ap1}) and (\ref{eq:BC_ap2}),
it comes
\begin{eqnarray*}
\varepsilon^2\delta  \left(\frac{\partial}{\partial t}\int_z^{\delta \eta} w\ dz + \delta\frac{\partial }{\partial x}\int_z^{\delta \eta} (uw)\ dz \right) - \varepsilon^2\delta^2 w^2 +
\delta p^a - p = -(\delta \eta - z)\nonumber\\
 + \varepsilon\delta\int_z^{\delta \eta} \frac{\partial}{\partial x}\left(\nu_0 \frac{\partial u}{\partial
z}\right)dz - 2\varepsilon\delta\nu_0\frac{\partial w}{\partial z}.
\label{eq:NS3_1}
\end{eqnarray*}
Classically we have
\begin{equation}
\frac{\partial u_s}{\partial x} = \left.\frac{\partial u}{\partial
x}\right|_s + \delta\frac{\partial \eta}{\partial
x}\left.\frac{\partial u}{\partial z}\right|_s = \left.\frac{\partial
u}{\partial x}\right|_s + {\cal O}(\varepsilon^2\delta),
\label{eq:deriv_comp}
\end{equation}
and using relations (\ref{eq:BC_ap3}), (\ref{eq:deriv_comp}) and the
Liebniz rule we have
\begin{equation*}
\varepsilon\delta\int_z^{\delta \eta} \frac{\partial}{\partial x}\left(\nu_0 \frac{\partial u}{\partial
z}\right) dz - 2\varepsilon\delta\nu_0\frac{\partial w}{\partial z} 
= \varepsilon\delta\nu_0\frac{\partial u}{\partial x} + \left.\varepsilon\delta\nu_0\frac{\partial u}{\partial x}\right|_s + {\cal O}(\varepsilon^3\delta).
\end{equation*}
This leads to the expression for the pressure $p$
\begin{eqnarray}
p = \delta p^a +(\delta \eta - z) +\varepsilon^2\delta  \left(\frac{\partial}{\partial t}\int_z^{\delta \eta} w\ dz + \delta\frac{\partial }{\partial x}\int_z^{\delta \eta} (uw)\ dz \right) \nonumber\\
- \varepsilon^2\delta^2 w^2 - \varepsilon\delta\nu_0\frac{\partial u}{\partial x} - \varepsilon\delta\nu_0\left.\frac{\partial u}{\partial x}\right|_s + {\cal O}(\varepsilon^3\delta).
\label{eq:pfin}
\end{eqnarray}
Hereafter several models of growing accuracy and complexity will be
derived, depending on the level of approximation chosen for
Eq.~(\ref{eq:pfin}). In the hydrostatic case, we will consider an
approximation of $p$ in ${\cal O}(\varepsilon^2\delta)$, then in
section~\ref{sec:SV_nhyd1} we will use two expressions of $p$
respectively in ${\cal O}(\varepsilon^2\delta^2,\varepsilon^3\delta)$ and in ${\cal O}(\varepsilon^3\delta)$.
\begin{remark}
For the derivation of Eq.~(\ref{eq:momentum}) note that due to the rescaling applied to the time derivative of $z_b$, we have
$$\int_{z_b}^{\delta\eta} \frac{\partial u}{\partial t} dz = \frac{\partial}{\partial t}\int_{z_b}^{\delta\eta} u\ dz - \delta\frac{\partial \eta}{\partial t}u_s + \varepsilon\delta\frac{\partial z_b}{\partial t}u_b.$$
\end{remark}
\begin{remark}
The second relation in (\ref{eq:BC1_ubis}) is crucial for the derivation of shallow
 water models. When considering large friction coefficients then the
 assumption of asymptotic regime $\hat{\kappa}=\varepsilon\kappa_0$ no
more holds and relation (\ref{eq:BC_ap5}) leads to
$$\frac{\partial u}{\partial z} = {\cal O}(1),$$
meaning the assumption of motion by slices has to be justified by other arguments. 
\end{remark}

\subsection{Hydrostatic approximation}
\label{subsec:SV}

We begin with the classical hydrostatic approximation. The objectives
of this paragraph are twofold. First we want to obtain the expression of
$\bar{u}$ as a function of $\delta$, $\varepsilon$, $\nu_0$, $\kappa_0$
and $H_\delta$. And second, we aim at verifying that despite the
parameter $\delta$, we recover the well-known formulation of the viscous
Saint-Venant system with friction as expressed in the following proposition
\begin{proposition}
The viscous Saint-Venant system defined by
\begin{eqnarray}
\lefteqn{\frac{\partial H}{\partial t} + \frac{\partial}{\partial x} \bigl(H\bar{u}\bigr) = 0,\label{eq:final11_1}}\\
\lefteqn{\frac{\partial (H\bar{u})}{\partial t} + \frac{\partial (H\bar{u}^2)}{\partial x} + \frac{g}{2}\frac{\partial H^2}{\partial x} = -H \frac{\partial p^a}{\partial x} -gH \frac{\partial z_b}{\partial x} + \frac{\partial}{\partial x}\bigl(4 \nu H\frac{\partial \bar{u}}{\partial x}\bigr)\nonumber}\\
& & \hspace*{4.7cm} - \frac{\kappa(\bar{{\bf v}},H)}{1 + \frac{\kappa(\bar{{\bf v}},H)}{3\nu}H} \bar{u},\hspace*{3cm}\label{eq:final11_2}
\end{eqnarray}
where $H = \eta - z_b$ and $\bar{{\bf v}}=(1,\frac{\partial z_b}{\partial x})^T\bar{u}$, results from an hydrostatic approximation in ${\cal O}(\varepsilon^2\delta)$ of the
 Navier-Stokes equations.
\label{prop:SVhydro}
\end{proposition}
{\it Proof of prop.~\ref{prop:SVhydro}}: we retain only the terms up to $\varepsilon\delta$ in the
expression (\ref{eq:pfin}) for the pressure $p$ i.e. we have
\begin{equation}
p = \delta p^a + (\delta\eta - z) - \varepsilon\delta\nu_0\frac{\partial u}{\partial
x} - \varepsilon\delta\nu_0 \left.\frac{\partial u}{\partial x}\right|_s + {\cal
O}(\varepsilon^2\delta).
\label{eq:phyd}
\end{equation}
And Eq.~(\ref{eq:momentum}) with Eqs.~(\ref{eq:u2}) and (\ref{eq:phyd}) gives
\begin{eqnarray*}
\lefteqn{\varepsilon\delta\frac{\partial \bigl(H_\delta\bar{u}\bigr)}{\partial t} + \varepsilon\delta^2\frac{\partial \bigl(H_\delta \overline{u^2}\bigr)}{\partial x} + \frac{\varepsilon}{2}\frac{\partial H_\delta^2}{\partial x} =}\nonumber\\ 
& & - \varepsilon\delta\kappa_0 u_b - \varepsilon\delta\frac{\partial }{\partial x}(H_\delta p^a) + \varepsilon\delta^2\frac{\partial \eta}{\partial x}p^a - \varepsilon\frac{\partial z_b}{\partial x} p_b + {\cal O}(\varepsilon^2\delta),
\end{eqnarray*}
that is also using the expression of $p$ obtained in Eq.~(\ref{eq:phyd})
\begin{eqnarray}
\lefteqn{\delta\frac{\partial \bigl(H_\delta\bar{u}\bigr)}{\partial t} + \delta^2\frac{\partial \bigl(H_\delta \overline{u^2}\bigr)}{\partial x} + \frac{1}{2}\frac{\partial H_\delta^2}{\partial x} =}\nonumber\\ 
& & - \delta\kappa_0 u_b - \delta H_\delta\frac{\partial p^a}{\partial x} - \frac{\partial z_b}{\partial x} p_b + {\cal O}(\varepsilon\delta).
\label{eq:momentumbis}
\end{eqnarray}
Note that due to the assumption concerning the time derivative of $z_b$
 and the associated rescaling, the first term in the left hand side of
 (\ref{eq:momentumbis}) reads 
$$\frac{\partial \bigl(H_\delta\bar{u}\bigr)}{\partial t} = H_\delta\frac{\partial \bar{u}}{\partial t} + \delta\frac{\partial (\eta - z_b)}{\partial t}\bar{u},$$
and (\ref{eq:momentumbis}) coupled with (\ref{eq:conserv}) gives
$$\delta\frac{\partial \bar{u}}{\partial t} +
\delta^2\bar{u}\frac{\partial \bar{u}}{\partial x} +
\delta\frac{\partial \eta}{\partial x} = -
\frac{\delta\kappa_0}{H_\delta} u_b - \delta\frac{\partial p^a}{\partial x} + {\cal O}(\varepsilon\delta).$$
Now we come back to Eq.~(\ref{eq:NS_2d_new2}), using
 (\ref{eq:u2}), (\ref{eq:phyd}) and (\ref{eq:momentumbis}) we get
\begin{eqnarray}
\delta\frac{\partial }{\partial z}\bigl(\nu_0\frac{\partial u}{\partial z}\bigr) & = & \varepsilon\delta\frac{\partial u}{\partial t} + \varepsilon\delta^2 u \frac{\partial u}{\partial x} + \varepsilon\delta^2 w \frac{\partial u}{\partial z} + \varepsilon\frac{\partial p}{\partial x} - \varepsilon^2\delta\frac{\partial }{\partial x}\left(\nu_0\frac{\partial u}{\partial x}\right)\nonumber\\
& = & \varepsilon\delta\frac{\partial \bar{u}}{\partial t} + \varepsilon\delta^2 \bar{u} \frac{\partial \bar{u}}{\partial x} + \varepsilon\delta\frac{\partial }{\partial x}(\eta + p^a) + {\cal O}(\varepsilon^2\delta)\nonumber \\
& = & -\frac{\varepsilon\delta\kappa_0}{H_\delta} u_b + {\cal O}(\varepsilon^2\delta). \label{eq:u_est}
\end{eqnarray}
Integrating from $z_b$ to $z$ and taking into account the boundary
condition (\ref{eq:BC_ap5}), we deduce
\begin{equation}
\frac{\partial u}{\partial z} = \frac{\varepsilon\kappa_0}{\nu_0} \left( 1 -
\frac{z -z_b}{H_\delta}\right) u_b + {\cal O}(\varepsilon^2),
\label{eq:duzb}
\end{equation}
and we obtain the following formula which gives an expression of the
vertical velocity though a parabolic correction
\begin{equation}
u = \left(1 + \frac{\varepsilon\kappa_0}{\nu_0} \bigl( z - z_b -
\frac{(z -z_b)^2}{2H_\delta}\bigr)\right) u_b + {\cal O}(\varepsilon^2).
\label{eq:u_approx}
\end{equation}
Then integrating from $z_b$ to $\delta \eta$, we obtain
\begin{equation}
\bar{u} = \left(1 + \frac{\varepsilon\kappa_0}{3\nu_0}H_\delta\right) u_b + {\cal O}(\varepsilon^2).
\label{eq:u_bar}
\end{equation}
Moreover
$$u^2 = \left(1 + \frac{2\varepsilon\kappa_0}{\nu_0} \bigl( z - z_b -
\frac{(z -z_b)^2}{2H_\delta}\bigr)\right)u_b^2 + {\cal O}(\varepsilon^2),$$
which yields
$$\overline{u^2} = \left(1 + \frac{2\varepsilon\kappa_0}{3\nu_0}H_\delta\right)u_b^2 + {\cal O}(\varepsilon^2),$$
meaning
\begin{equation}
\overline{u^2} = \bar{u}^2 + {\cal O}(\varepsilon^2).
\label{eq:u2_bis}
\end{equation}
Using (\ref{eq:phyd}), (\ref{eq:u_approx}) and (\ref{eq:u_bar}), the right hand side of Eq.~(\ref{eq:momentum}) can be written
\begin{eqnarray}
\varepsilon\delta\frac{\partial }{\partial x}\int_{z_b}^{\delta \eta} 2\nu_0\frac{\partial u}{\partial x} dz + \delta^2 \frac{\partial \eta}{\partial x} p^a - \frac{\partial z_b}{\partial x} \left(p_b + 2\varepsilon\delta\nu_0\left.\frac{\partial u}{\partial x}\right|_b\right)\nonumber\\
 - \delta\kappa_0 \left(1 + \frac{5\varepsilon^2}{2}\left(\frac{\partial z_b}{\partial x}\right)^2\right) u_b = -\delta\kappa_0 u_b - H_\delta\frac{\partial z_b}{\partial x}\nonumber\\
 + \delta \frac{\partial H_\delta}{\partial x}p^a + \varepsilon\delta\frac{\partial}{\partial x}\bigl(2 \nu_0 H_\delta\frac{\partial \bar{u}}{\partial x}\bigr) + {\cal O}(\varepsilon^2\delta).\label{eq:BC_b}
\end{eqnarray}
Finally from Eqs.~(\ref{eq:conserv}), (\ref{eq:momentum}), (\ref{eq:u_bar}), (\ref{eq:u2_bis}) and (\ref{eq:BC_b}), we obtain the model
\begin{eqnarray*}
\frac{\partial \eta}{\partial t} - \frac{\partial z_b}{\partial t} + \frac{\partial}{\partial x} \bigl(H_\delta\bar{u}\bigr) = 0,\qquad\label{eq:final1_1}\\
\delta\frac{\partial (H_\delta\bar{u})}{\partial t} + \delta^2\frac{\partial (H_\delta \bar{u}^2)}{\partial x} + \frac{1}{2}\frac{\partial H_\delta^2}{\partial x} & = & -H_\delta \frac{\partial }{\partial x}(z_b + \delta p^a) - \frac{\delta\kappa_0}{1 + \frac{\varepsilon\kappa_0}{3\nu_0}H_\delta} \bar{u} \nonumber\\
& & + \varepsilon\delta\frac{\partial}{\partial x}\bigl(4 \nu_0 H_\delta\frac{\partial \bar{u}}{\partial x}\bigr) + {\cal O}(\varepsilon^2\delta). \label{eq:final1_2}
\end{eqnarray*}
In terms of the initial variables, the preceding model becomes
(\ref{eq:final11_1})-(\ref{eq:final11_2}) that complete the proof of prop.~\ref{prop:SVhydro}. Note that when the bathymetry
is constant $z_b(x,t)=z_b^0$, this formulation is
equivalent to the viscous Saint-Venant system obtained by Gerbeau {\it
et al.} \cite{gerbeau} and Ferrari {\it et al.} \cite{saleri}.

\section{Two non-hydrostatic shallow water models}
\label{sec:SV_nhyd1}

In the previous paragraph we have obtained an approximation of the
Navier-Stokes equations up to $\varepsilon\delta$ terms using an
hydrostatic approximation of the pressure $p$. In this section
we consider two more acurate approximations of the pressure $p$
respectively in ${\cal O}(\varepsilon^2\delta^2)$ and ${\cal
O}(\varepsilon^3\delta)$ leading to two non-hydrostatic extensions of
the Saint-Venant system.

\subsection{First extension, $\delta \ll 1$}
\label{subsec:extension_1}

The first refinement of the classical Saint-Venant model (\ref{eq:final11_1})-(\ref{eq:final11_2}) is achieved by considering the pressure
$p$ given by Eq.~(\ref{eq:pfin}) with the terms up to ${\cal
O}(\varepsilon^2\delta^2)$. This means we consider the momentum equation
along $z$ is no more reduced to
$$\frac{\partial p}{\partial z} = -1 + \frac{\partial }{\partial x}\left( \varepsilon\delta\nu_0\frac{\partial u}{\partial z} \right) + \varepsilon\delta\frac{\partial }{\partial z}\left(2\nu_0\frac{\partial w}{\partial z}\right) + {\cal O}(\varepsilon^2\delta),$$
but given by
$$\varepsilon^2\delta\frac{\partial w}{\partial t} + \frac{\partial p}{\partial z} = -1 + \frac{\partial }{\partial x}\left( \varepsilon\delta\nu_0\frac{\partial u}{\partial z} \right) + \varepsilon\delta\frac{\partial }{\partial z}\left(2\nu_0\frac{\partial w}{\partial z}\right) + {\cal O}(\varepsilon^2\delta^2),$$
and the convective terms are still neglected. Since we keep the terms in
$\varepsilon^2\delta$ and drop those in $\varepsilon^2\delta^2$, this
means we assume $\delta \ll 1$ and due to the applied rescaling this implies $U \ll C$ so we are in a fluvial regime. The following result holds.
\begin{proposition}
The system defined by
\begin{eqnarray}
\lefteqn{\frac{\partial H}{\partial t} + \frac{\partial}{\partial x} \bigl(H\bar{u}\bigr) = 0,} \label{eq:sv_nonhyd_init1}\\
\lefteqn{\frac{\partial}{\partial t}(H \bar{u}) + \frac{\partial}{\partial x} (H \bar{u}^2) + \frac{\partial}{\partial x}\left(\frac{g}{2}H^2 -\frac{z_b^3}{6}\frac{\partial^2 \bar{u}}{\partial x \partial t} + \frac{z_b^2}{2}\frac{\partial^2 (z_b\bar{u})}{\partial x\partial t} \right) = \nonumber}\\
& & - H\frac{\partial p^a}{\partial x} + \frac{\partial }{\partial x}\left(4\nu H \frac{\partial \bar{u}}{\partial x} + \frac{\kappa(\bar{{\bf v}},H)}{6}z_b\left(z_b\frac{\partial \bar{u}}{\partial x} + 7\frac{\partial z_b}{\partial x}\bar{u}\right)\right) \nonumber\\
& & - \frac{\kappa(\bar{{\bf v}},H)}{2}\frac{\partial z_b}{\partial x}\left(z_b\frac{\partial \bar{u}}{\partial x}-\frac{\partial z_b}{\partial x}\bar{u}\right) - \frac{\partial z_b}{\partial x}\left( gH +\frac{z_b^2}{2}\frac{\partial^2 \bar{u}}{\partial x\partial t} - z_b\frac{\partial^2 (z_b\bar{u})}{\partial x \partial t}\right)\nonumber\\
& & -\frac{\kappa(\bar{{\bf v}},H)}{1+\frac{\kappa(\bar{{\bf v}},H) H}{3\nu}}\left(1 + \frac{5}{2}\left(\frac{\partial z_b}{\partial x}\right)^2\right)\bar{u} - \frac{z_b^2}{2}\frac{\partial^3 z_b}{\partial x\partial t^2},\label{eq:sv_nonhyd_init2}
\end{eqnarray}
where $\bar{{\bf v}}=(1,\frac{\partial z_b}{\partial x})^T\frac{\bar{u}}{1+\frac{\kappa_l H}{3\nu}}$ results from an approximation in ${\cal O}(\varepsilon^2\delta^2,\varepsilon^3\delta)$ of the Navier-Stokes equations.
\label{prop:SV+}
\end{proposition}

The proof of proposition~\ref{prop:SV+} is given in the next paragraph,
we examine here some properties of the model (\ref{eq:sv_nonhyd_init1})-(\ref{eq:sv_nonhyd_init2}).

Note that except for the dissipative terms corresponding to viscosity or friction,
all the terms added in the non-hydrostatic model
(\ref{eq:sv_nonhyd_init1})-(\ref{eq:sv_nonhyd_init2}) compared to the original Saint-Venant model
(\ref{eq:final11_1})-(\ref{eq:final11_2}) appear
as time derivative of the variables $z_b$, $\eta$ or $\bar{u}$. This means in
a stationary regime, the solutions of (\ref{eq:final11_1})-(\ref{eq:final11_2}) and
(\ref{eq:sv_nonhyd_init1})-(\ref{eq:sv_nonhyd_init2}) are identical.

We first examine the system
(\ref{eq:sv_nonhyd_init1})-(\ref{eq:sv_nonhyd_init2}) without friction
and viscosity. Starting from the Euler equations instead of the
Navier-Stokes equations does not allow to account for the motion by slices
as obtained in relations (\ref{eq:BC1_uter}) and (\ref{eq:duzb}). So if
one wants to neglect the viscosity and friction effects in the model
(\ref{eq:sv_nonhyd_init1})-(\ref{eq:sv_nonhyd_init2}), it is necessary
to consider an asymptotic regime for example under the form $\nu =
\beta\nu_{nv}$, $\kappa = \beta^2 \kappa_{nf}$ --~and conversely $\nu_0
= \beta\nu_{0,nv}$, $\kappa_0 = \beta^2 \kappa_{0,nf}$~-- with $\beta \ll 1$.
Introducing the preceding asymptotic regime and considering $\beta \rightarrow 0$, the formulation of (\ref{eq:sv_nonhyd_init1})-(\ref{eq:sv_nonhyd_init2}) reads
\begin{eqnarray*}
\lefteqn{\frac{\partial H}{\partial t} + \frac{\partial}{\partial x} \bigl(H\bar{u}\bigr) = 0,}\\
\lefteqn{\frac{\partial}{\partial t}(H \bar{u}) + \frac{\partial}{\partial x} (H \bar{u}^2) + \frac{\partial}{\partial x}\left(\frac{g}{2}H^2 -\frac{z_b^3}{6}\frac{\partial^2 \bar{u}}{\partial x \partial t} + \frac{z_b^2}{2}\frac{\partial^2 (z_b\bar{u})}{\partial x\partial t} \right) = \nonumber}\\
& & \frac{\partial z_b}{\partial x}\left( -gH -\frac{z_b^2}{2}\frac{\partial^2 \bar{u}}{\partial x\partial t} + z_b\frac{\partial^2 (z_b\bar{u})}{\partial x \partial t}\right) - H\frac{\partial p^a}{\partial x} - \frac{z_b^2}{2}\frac{\partial^3 z_b}{\partial x\partial t^2}.
\end{eqnarray*}
or equivalently in a non-conservative form
\begin{eqnarray*}
& & \frac{\partial H}{\partial t} + \frac{\partial}{\partial x} \bigl(H\bar{u}\bigr) = 0,\label{eq:sv_nonhyd_noncons1}\\
&& \frac{\partial \bar{u}}{\partial t} + \bar{u}\frac{\partial \bar{u}}{\partial x} + g\frac{\partial \eta}{\partial x} + \frac{z_b^2}{6}\frac{\partial^3 \bar{u}}{\partial x^2 \partial t} - \frac{z_b}{2}\frac{\partial^3 (z_b\bar{u})}{\partial x^2\partial t} = - \frac{\partial p^a}{\partial x} + \frac{z_b}{2}\frac{\partial^3 z_b}{\partial x\partial t^2},\label{eq:sv_nonhyd_noncons2}
\end{eqnarray*}
that is analogous to the expression obtained by Peregrine
\cite{peregrine}. It is worth being noticed that, in any case, the formulations
obtained by Nwogu \cite{nwogu}, Walkley \cite{walkley}, Saut {\it et
al.} \cite{saut1} and Soares Frazao {\it et al.} \cite{zech} are different from the preceding ones. The differences
lie either in the continuity equation or in the momentum equation.

The mathematical and numerical analysis of the obtained model is
not in the scope of this paper but let us mention some interesting works
in the literature. The Sobolev equation
$$-\frac{\partial}{\partial x}(a(x)\frac{\partial^2 u }{\partial x \partial
t}) + c(x)\frac{\partial u}{\partial t} = -\frac{\partial}{\partial
x}(\alpha(x) \frac{\partial u}{\partial x}) + \beta(x) \frac{\partial u}{\partial x},$$
has been studied by several authors \cite{arnold,bbm} as an alternative
to the Korteweg-de Vries equations. Perotto and Saleri \cite{perotto}
proposed an {\it a posteriori} error analysis for the Peregrine
formulation of the Boussinesq system with constant bathymetry.
Bona {\it et al.} \cite{saut1,saut2} have studied the well-posedness of
several high-order generalizations of the Boussinesq equations.

\subsection{Derivation}
\label{subsec:nhyd1}

{\it Proof of prop.~\ref{prop:SV+}:} the refinement of the classical Saint-Venant model (\ref{eq:final11_1})-(\ref{eq:final11_2}) is achieved by improving the approximation for the pressure
$p$. Actually, if we only drop the terms in ${\cal O}(\varepsilon^2\delta^2)$ in the
momentum equation along $z$ so the system (\ref{eq:NS_2d_re1})-(\ref{eq:NS_2d_re3}) becomes
\begin{eqnarray*}
&& w = \frac{\partial z_b}{\partial t} - \frac{\partial}{\partial x} \int_{z_b}^z u\ dz,\label{eq:NS_2d_re_new1}\\
&& \varepsilon\delta\frac{\partial u}{\partial t} + \varepsilon\delta^2 \frac{\partial u^2}{\partial x} + \varepsilon\delta^2 \frac{\partial u w}{\partial z} + \varepsilon\frac{\partial p}{\partial x} = \varepsilon^2\delta\frac{\partial }{\partial x}\left(2\nu_0\frac{\partial u}{\partial x}\right) + \frac{\partial }{\partial z}\left( \delta\nu_0\frac{\partial u}{\partial z}\right),\label{eq:NS_2d_re_new2}\\
&& \varepsilon^2\delta\frac{\partial w}{\partial t} + \frac{\partial p}{\partial z} = -1 + \frac{\partial }{\partial x}\left( \varepsilon\delta\nu_0\frac{\partial u}{\partial z} \right) + \varepsilon\delta\frac{\partial }{\partial z}\left(2\nu_0\frac{\partial w}{\partial z}\right) + {\cal O}(\varepsilon^2\delta^2),\label{eq:NS_2d_re_new3}
\end{eqnarray*}
with the boundary conditions (\ref{eq:BC_ap1})-(\ref{eq:BC_ap5}). This means we consider the pressure $p$ is given by (\ref{eq:pfin}) where we retain
only the terms up to $\varepsilon^2\delta^2$ and $\varepsilon^3\delta$ i.e.
\begin{eqnarray*}
p_{nh} & = & \delta p^a + (\delta\eta - z) - \varepsilon\delta\nu_0\left.\frac{\partial u}{\partial x}\right|_s - \varepsilon\delta\nu_0 \frac{\partial u}{\partial x} + \varepsilon^2\delta\frac{\partial }{\partial t}\int_z^{\delta \eta} w\ dz \\
&& + {\cal O}(\varepsilon^2\delta^2,\varepsilon^3\delta),
\end{eqnarray*}
leading to
\begin{eqnarray}
p_{nh} & = & \delta p^a + (\delta\eta - z) - \varepsilon\delta\nu_0\left.\frac{\partial u}{\partial
x}\right|_s - \varepsilon\delta\nu_0 \frac{\partial u}{\partial
x} + \varepsilon^2\delta (\delta \eta - z) \frac{\partial^2 z_b}{\partial t^2}\nonumber\\
& & - \varepsilon^2\delta\frac{\partial}{\partial t}\int_z^{\delta \eta} \frac{\partial}{\partial x} \int_{z_b}^{z} u dz_1 dz + {\cal O}(\varepsilon^2\delta^2,\varepsilon^3\delta).
\label{eq:pres}
\end{eqnarray}
Retaining only the terms up to ${\cal
O}(\varepsilon^2\delta^2,\varepsilon^3\delta)$, relation (\ref{eq:momentum}) gives
\begin{eqnarray}
\lefteqn{\delta\frac{\partial}{\partial t}\int_{z_b}^{\delta \eta} u\ dz + \delta^2\frac{\partial}{\partial x}\int_{z_b}^{\delta \eta} u^2\ dz + \frac{\partial}{\partial x}\int_{z_b}^{\delta \eta} p_{nh}\ dz }\nonumber\\
& & = \varepsilon\delta\frac{\partial }{\partial x}\int_{z_b}^{\delta \eta} 2\nu_0\frac{\partial u}{\partial x} dz + \delta^2 \frac{\partial \eta}{\partial x} p^a + \frac{\partial z_b}{\partial x} \left(2\varepsilon\delta\nu_0\left.\frac{\partial w}{\partial z}\right|_b - \left. p_{nh}\right|_b \right) \nonumber\\
& &\quad - \delta\kappa_0 \left(1 + \frac{5\varepsilon^2}{2}\bigl(\frac{\partial z_b}{\partial x}\bigr)^2\right) u_b + {\cal O}(\varepsilon^2\delta^2,\varepsilon^3\delta).
\label{eq:momentum1}
\end{eqnarray}
Now we derive the expressions for the quantities appearing in
 (\ref{eq:pres}) and (\ref{eq:momentum1}) and depending on $u$, $w$ and
 $p$. Since $\kappa_0 = \kappa_{0,l} + {\cal O}(\varepsilon)$, from Eqs.~(\ref{eq:u_approx}) and (\ref{eq:u_bar}) we have
\begin{eqnarray*}
\frac{\partial u}{\partial x} & = & \left(1 + \frac{\varepsilon\kappa_0}{\nu_0} \bigl( z - z_b -
\frac{(z -z_b)^2}{2H_\delta}\bigr)\right) \frac{\partial u_b}{\partial x}\nonumber\\
& & + \frac{\varepsilon\kappa_0}{\nu_0} \left( \frac{\partial z_b}{\partial x}\bigl(-1 + \frac{z -z_b}{H_\delta}\bigr) +
\frac{\partial H_\delta}{\partial x}\frac{(z -z_b)^2}{2H_\delta^2}\right) u_b + {\cal O}(\varepsilon^2)\nonumber\\
& = & \left(1 + \frac{\varepsilon\kappa_0}{\nu_0} \bigl( z - z_b - \frac{(z -z_b)^2}{2H_\delta}\bigr)\right) \frac{\partial u_b}{\partial x}\nonumber\\
& & - \frac{\varepsilon\kappa_0}{\nu_0} \frac{\partial z_b}{\partial x}\left(1 - \frac{z -z_b}{H_\delta} + \frac{(z -z_b)^2}{2H_\delta^2}\right) u_b + {\cal O}(\varepsilon\delta),
\end{eqnarray*}
so
\begin{eqnarray}
\left.\frac{\partial u}{\partial x}\right|_s & = & \bigl(1 + \frac{\varepsilon\kappa_0}{2\nu_0} H_\delta\bigr)\frac{\partial u_b}{\partial x} + \frac{\varepsilon\kappa_0}{2\nu_0}\frac{\partial H_\delta}{\partial x}u_b + {\cal O}(\varepsilon^2)\label{eq:dus0}\\
& = & \bigl(1 + \frac{\varepsilon\kappa_0}{6\nu_0} H_\delta\bigr)\frac{\partial \bar{u}}{\partial x} + \frac{\varepsilon\kappa_0}{6\nu_0}\frac{\partial H_\delta}{\partial x}\bar{u} + {\cal O}(\varepsilon^2),\label{eq:dus}\\
\left.\frac{\partial u}{\partial x}\right|_b & = & \frac{\partial u_b}{\partial x} - \frac{\varepsilon\kappa_0}{\nu_0}\frac{\partial z_b}{\partial x}u_b + {\cal O}(\varepsilon^2)\label{eq:dub0}\\
& = & \bigl(1 - \frac{\varepsilon\kappa_0}{3\nu_0} H_\delta\bigr)\frac{\partial \bar{u}}{\partial x} - \frac{\varepsilon\kappa_0}{\nu_0}\bigl(\frac{\partial z_b}{\partial x} + \frac{1}{3}\frac{\partial H_\delta}{\partial x}\bigr)\bar{u} + {\cal O}(\varepsilon^2),\label{eq:dub}
\end{eqnarray}
and
\begin{eqnarray}
\int_{z_b}^{\delta\eta} \nu_0 \frac{\partial u}{\partial x}  & = & \nu_0 H_\delta \bigl( 1 + \frac{\varepsilon\kappa_0}{3 \nu_0}H_\delta \bigr) \frac{\partial u_b}{\partial x} - \frac{\varepsilon\kappa_0 H_\delta}{2}\bigl(\frac{\partial z_b}{\partial x} - \frac{1}{3}\frac{\partial H_\delta}{\partial x}\bigr)u_b + {\cal O}(\varepsilon^2)\nonumber\\
& = & \nu_0 H_\delta\frac{\partial \bar{u}}{\partial x} - \frac{\varepsilon\kappa_0}{3}H_\delta\frac{\partial H_\delta}{\partial x}\bar{u} - \frac{\varepsilon\kappa_0 H_\delta}{2}\bigl(\frac{\partial z_b}{\partial x} - \frac{1}{6}\frac{\partial H_\delta}{\partial x}\bigr)\bar{u} + {\cal O}(\varepsilon^2) \nonumber\\
& = & \nu_0 H_\delta\frac{\partial \bar{u}}{\partial x} - \frac{\varepsilon\kappa_0}{2} H_\delta\left(\frac{\partial z_b}{\partial x} + \frac{1}{3}\frac{\partial H_\delta}{\partial x}\right)\bar{u} + {\cal O}(\varepsilon^2), \label{eq:du}
\end{eqnarray}
and finally from (\ref{eq:deriv_comp}) we get
\begin{eqnarray*}
\int_{z_b}^{\delta\eta} \nu_0 \left.\frac{\partial u}{\partial x}\right|_s & = & \int_{z_b}^{\delta\eta} \nu_0 \frac{\partial u_s}{\partial x} + {\cal O}(\varepsilon^2\delta)\nonumber\\
& = & \nu_0 H_\delta \left( 1 + \frac{\varepsilon\kappa_0}{6 \nu_0}H_\delta \right) \frac{\partial \bar{u}}{\partial x} + \frac{\varepsilon\kappa_0}{6} \frac{\partial H_\delta}{\partial x}H_\delta \bar{u} + {\cal O}(\varepsilon^2).
\label{eq:dusint}
\end{eqnarray*}
From (\ref{eq:dus}) and (\ref{eq:dub}) we have
\begin{eqnarray}
\left. p_h\right|_b - 2\varepsilon\delta\nu_0\left. \frac{\partial w}{\partial z}\right|_b & = & \delta p^a + H_\delta + \varepsilon\delta\nu_0\left. \frac{\partial u}{\partial x}\right|_b - \varepsilon\delta\nu_0\left.\frac{\partial u}{\partial
x}\right|_s \nonumber\\
& = & \delta p^a + H_\delta - \varepsilon^2\delta\frac{\kappa_0}{2}H_\delta\frac{\partial \bar{u}}{\partial x}\nonumber\\
& & - \varepsilon^2\delta\kappa_0\left( \frac{\partial z_b}{\partial x}+\frac{1}{2}\frac{\partial H_\delta}{\partial x}\right)\bar{u} + {\cal O}(\varepsilon^3\delta), \label{eq:pb}
\end{eqnarray}
and
\begin{eqnarray}
\int_{z_b}^{\delta\eta} \bigl(2\varepsilon\delta\nu_0\frac{\partial u}{\partial x} - p_h \bigr) dz  =  -H_\delta \delta p^a - \frac{H_\delta^2}{2} +4\varepsilon\delta\nu_0 H_\delta \frac{\partial \bar{u}}{\partial x}\nonumber\\
+ \varepsilon^2\delta\kappa_0 H_\delta\left( \frac{H_\delta}{6}\frac{\partial \bar{u}}{\partial x} - \frac{7}{6}\frac{\partial z_b}{\partial x}\bar{u} - \frac{\delta}{3}\frac{\partial \eta}{\partial x}\bar{u}\right) + {\cal O}(\varepsilon^3\delta),
\label{eq:pm}
\end{eqnarray}
where $p_h$ corresponds to the gravitational, viscous and friction part of the pressure $p$ given by Eq.~(\ref{eq:pres}) i.e.
\begin{equation}
p_h = \delta p^a + (\delta\eta - z) - \varepsilon\delta\nu_0\left.\frac{\partial u}{\partial
x}\right|_s - \varepsilon\delta\nu_0 \frac{\partial u}{\partial x}.
\label{eq:ph}
\end{equation}
Inserting (\ref{eq:du}), (\ref{eq:pb}) and (\ref{eq:pm}) in equilibrium
(\ref{eq:momentum1}) leads to
\begin{eqnarray*}
\lefteqn{\delta\frac{\partial}{\partial t}(H_\delta \bar{u}) + \delta^2\frac{\partial}{\partial x} (H_\delta \bar{u}^2) + \frac{1}{2}\frac{\partial H_\delta^2}{\partial x} + \frac{\partial }{\partial x}\int_{z_b}^{\delta \eta} \Delta p_{nh}\ dz = }\nonumber\\
& & - H_\delta\frac{\partial }{\partial x}(\delta p^a + z_b) + \varepsilon\delta\frac{\partial }{\partial x}\bigl(4\nu_0H_\delta \frac{\partial \bar{u}}{\partial x}\bigr) \nonumber\\
& & + \frac{\varepsilon^2\delta}{6}\frac{\partial}{\partial x}\left( \kappa_0 z_b\left(z_b\frac{\partial \bar{u}}{\partial x} + 7\frac{\partial z_b}{\partial x}\bar{u}\right)\right) - \frac{\varepsilon^2\delta\kappa_0}{2}\frac{\partial z_b}{\partial x}\left(-\frac{\partial z_b}{\partial x}\bar{u} + z_b\frac{\partial \bar{u}}{\partial x}\right)\nonumber\\
& & -\frac{\partial z_b}{\partial x}\left. \Delta p_{nh}\right|_b -\delta\kappa_0\bigl(1 + \frac{5\varepsilon^2}{2}\bigl(\frac{\partial z_b}{\partial x}\bigr)^2\bigr)u_b + {\cal O}(\varepsilon^2\delta^2,\varepsilon^3\delta), \label{eq:momentum3}
\end{eqnarray*}
where $\Delta p_{nh} = p_{nh} - p_h$. And using the expression for the pressure $p_{nh}$ given in
Eq.~(\ref{eq:pres}) it comes
\begin{eqnarray}
\int_{z_b}^{\delta\eta} \Delta p_{nh}\ dz & = & -\varepsilon^2\delta\frac{H_\delta^2}{6}(2\delta\eta + z_b)\frac{\partial^2 \bar{u}}{\partial x \partial t} + \varepsilon^2\delta\frac{H_\delta^2}{2}\frac{\partial^2 (z_b\bar{u})}{\partial x \partial t} \nonumber\\
& & - \varepsilon^2\delta^2 H_\delta \frac{\partial \eta}{\partial t} \left(\delta\eta\frac{\partial \bar{u}}{\partial x} - \frac{\partial (z_b\bar{u})}{\partial x} \right) + \varepsilon^2\delta\frac{H_\delta^2}{2}\frac{\partial^2 z_b}{\partial t^2}\label{eq:presint}\\
& = & \varepsilon^2\delta\frac{z_b^2}{2}\left(-\frac{z_b}{3}\frac{\partial^2 \bar{u}}{\partial x \partial t} + \frac{\partial^2 (z_b\bar{u})}{\partial x\partial t} + \frac{\partial^2 z_b}{\partial t^2}\right) + {\cal O}(\varepsilon^2\delta^2,\varepsilon^3\delta),\nonumber
\end{eqnarray}
and
\begin{eqnarray}
\left. \Delta p_{nh}\right|_b & = & -\varepsilon^2\delta\frac{\delta^2\eta^2-z_b^2}{2}\frac{\partial^2 \bar{u}}{\partial x\partial t} + \varepsilon^2\delta H_\delta\frac{\partial^2 (z_b\bar{u})}{\partial x\partial t} \nonumber\\
& & - \varepsilon^2\delta^2\frac{\partial \eta}{\partial t}\left(\delta\eta\frac{\partial \bar{u}}{\partial x} - \frac{\partial (z_b\bar{u})}{\partial x}\right) + \varepsilon^2\delta H_\delta\frac{\partial^2 z_b}{\partial t^2}\label{eq:presint2}\\
& = & \frac{\varepsilon^2\delta}{2}\left( z_b^2\frac{\partial^2 \bar{u}}{\partial x\partial t} - 2z_b\frac{\partial^2 (z_b\bar{u})}{\partial x\partial t}\right) - \varepsilon^2\delta z_b\frac{\partial^2 z_b}{\partial t^2} + {\cal O}(\varepsilon^2\delta^2,\varepsilon^3\delta).\nonumber
\end{eqnarray}
We finally obtain the model
\begin{eqnarray}
\lefteqn{\frac{\partial H_\delta}{\partial t} + \frac{\partial}{\partial x} \bigl(H_\delta\bar{u}\bigr) = 0,}\label{eq:sv_nonhyd1}\\
\lefteqn{\delta\frac{\partial}{\partial t}(H_\delta \bar{u}) + \delta^2\frac{\partial}{\partial x} (H_\delta \bar{u}^2) + \frac{1}{2}\frac{\partial H_\delta^2}{\partial x} - \varepsilon^2\delta\frac{\partial}{\partial x}\left(\frac{z_b^3}{6}\frac{\partial^2 \bar{u}}{\partial x \partial t} - \frac{z_b^2}{2}\frac{\partial^2 (z_b\bar{u})}{\partial x \partial t}\right) = }\nonumber\\
& & - H_\delta\frac{\partial }{\partial x}(\delta p^a + z_b) + \frac{\partial }{\partial x}\left(4\varepsilon\delta\nu_0H_\delta \frac{\partial \bar{u}}{\partial x} + \frac{\varepsilon\kappa_0}{6} z_b\left( z_b\frac{\partial \bar{u}}{\partial x} + 7\frac{\partial z_b}{\partial x}\bar{u}\right)\right) \nonumber\\
& & - \frac{\varepsilon^2\delta\kappa_0}{2}\frac{\partial z_b}{\partial x}\left(-\frac{\partial z_b}{\partial x}\bar{u} + z_b\frac{\partial \bar{u}}{\partial x}\right) + \varepsilon^2\delta\frac{\partial z_b}{\partial x}\left( -\frac{z_b^2}{2}\frac{\partial^2 \bar{u}}{\partial x\partial t} + z_b\frac{\partial^2 (z_b\bar{u})}{\partial x \partial t}\right) \nonumber\\
& & -\delta\kappa_0\left(1 + \frac{5\varepsilon^2}{2}\bigl(\frac{\partial z_b}{\partial x}\bigr)^2\right)u_b  - \varepsilon^2\delta\frac{z_b^2}{2}\frac{\partial^3 z_b}{\partial x\partial t^2} + {\cal O}(\varepsilon^2\delta^2,\varepsilon^3\delta),\label{eq:sv_nonhyd2}
\end{eqnarray}
that complete the proof of proposition~\ref{prop:SV+}. When the terms in ${\cal O}(\varepsilon^2\delta)$ are dropped in
(\ref{eq:sv_nonhyd2}), we verify that we recover the classical viscous hydrostatic
Saint-Venant model with friction (\ref{eq:final11_1})-(\ref{eq:final11_2}).

\subsection{Energy equality}
\label{subsec:prop}

Until now, we have not verified the derived models satisfy an energy
equality. The system (\ref{eq:final11_1})-(\ref{eq:final11_2}) that is equivalent to the Saint-Venant
system, admits a dissipation energy \cite{audusse,bouchut}. Indeed we have
\begin{eqnarray}
\frac{\partial E_h}{\partial t} + \frac{\partial}{\partial x}
\left(\bar{u}\bigl(E_h+g\frac{H^2}{2}\bigr) - 4\nu H\bar{u}\frac{\partial \bar{u}}{\partial x}\right) = -H\frac{\partial p^a}{\partial t} - 4\nu H\bigl(\frac{\partial \bar{u}}{\partial x}\bigr)^2\nonumber\\
 - \frac{\kappa(\bar{{\bf v}},H)}{1+\frac{\kappa(\bar{{\bf v}},H) H}{3 \nu}}\bar{u}^2  + gH\frac{\partial z_b}{\partial t},\label{eq:energy}
\end{eqnarray}
with $E_h=\frac{H\bar{u}^2}{2}+\frac{gH(\eta+z_b)}{2} + Hp^a$. The energy
equality (\ref{eq:energy}) associated with the hydrostatic Saint-Venant
model can be obtained using classical computations by multiplying Eq.~(\ref{eq:momentum}) when $p=p_h$ by the velocity $\bar{u}$.

The only differences between the hydrostatic Saint-Venant model (\ref{eq:final11_1})-(\ref{eq:final11_2}) and its extended version (\ref{eq:sv_nonhyd_init1})-(\ref{eq:sv_nonhyd_init2}) comes from
\begin{itemize}
\item[$\bullet$] the non hydrostatic terms of the pressure $p_{nh}$,
\item[$\bullet$] the terms involving the viscosity and the friction at the bottom,
\end{itemize}
so the energy equality for (\ref{eq:sv_nonhyd1})-(\ref{eq:sv_nonhyd2}) will differ from Eq.~(\ref{eq:energy}) only by the terms
\begin{eqnarray*}
{\cal C}_1 & = & \bar{u}\frac{\partial}{\partial x}\int_{z_b}^{\delta\eta}
\Delta p_{nh} +\bar{u}\frac{\partial z_b}{\partial x}
\left. \Delta p_{nh}\right|_b,\\
{\cal C}_2 & = & \bar{u}\frac{\partial}{\partial x}\int_{z_b}^{\delta\eta} \left(2\varepsilon\delta\nu_0\frac{\partial u}{\partial x}\right),\\
{\cal C}_3 & = & \bar{u}\frac{\partial}{\partial x}\int_{z_b}^{\delta\eta} p_{v,f},\\
{\cal C}_4 & = & \bar{u}\frac{\partial z_b}{\partial x} \left(2\varepsilon\delta\nu_0\left.\frac{\partial w}{\partial z}\right|_b - \left. p_{v,f}\right|_b - \varepsilon\delta\nu_0\frac{\partial z_b}{\partial x}\left.\frac{\partial u}{\partial z}\right|_b\right),
\end{eqnarray*}
where $\Delta p_{nh} = p_{nh} - p_h$ and $p_{v,f} = p_h - \delta p^a$
denotes the terms in the pressure $p$ containing the viscosity and
friction. The quantities ${\cal C}_1$-${\cal C}_4$ corresponding to the non-hydrostatic terms,
come from the multiplication of Eq.~(\ref{eq:momentum}) by
$\bar{u}$ and have to be added to (\ref{eq:energy}). Since $\bar{u}= u + {\cal
O}(\varepsilon)= u_b + {\cal O}(\varepsilon)$ and $\Delta p_{nh} = {\cal O}(\varepsilon^2\delta^2)$, we rewrite ${\cal C}_1$ under the form
\begin{eqnarray*}
{\cal C}_1 & = & \bar{u}\frac{\partial}{\partial x}\int_{z_b}^{\delta\eta} \Delta p_{nh} + u_b\frac{\partial z_b}{\partial x} \left. \Delta p_{nh}\right|_b + {\cal O}(\varepsilon^2\delta^2)\nonumber\\
& = & \frac{\partial}{\partial x}\int_{z_b}^{\delta\eta} u p_{nh} - \int_{z_b}^{\delta\eta} \frac{\partial u }{\partial x} \Delta p_{nh}  + u_b\frac{\partial z_b}{\partial x} \left. \Delta p_{nh}\right|_b + {\cal O}(\varepsilon^2\delta^2)\nonumber\\
& = & \frac{\partial}{\partial x} \int_{z_b}^{\delta\eta} u \Delta p_{nh} + [w \Delta p_{nh}]_{z_b}^{\delta\eta} - \int_{z_b}^{\delta\eta} w \frac{\partial \Delta p_{nh} }{\partial z}  + u_b\frac{\partial z_b}{\partial x} \left. \Delta p_{nh}\right|_b + {\cal O}(\varepsilon^2\delta^2)\nonumber\\
& = & \frac{\partial}{\partial x} \int_{z_b}^{\delta\eta} u \Delta p_{nh} + \left. w_s \Delta p_{nh}\right|_s - \int_{z_b}^{\delta\eta} w \frac{\partial \Delta p_{nh} }{\partial z}  - \frac{\partial z_b}{\partial t} \left. \Delta p_{nh}\right|_b + {\cal O}(\varepsilon^2\delta^2),\nonumber
\end{eqnarray*}
where relation (\ref{eq:BC_ap4})
has been used. From Eqs.~(\ref{eq:pres}) and (\ref{eq:ph}), we have
\begin{eqnarray*}
& & \left. \Delta p_{nh}\right|_s = {\cal O}(\varepsilon^2\delta^2),\\
& & \left. \Delta p_{nh}\right|_b = \varepsilon^2\delta\int_{z_b}^{\delta\eta} \frac{\partial w}{\partial t} + {\cal O}(\varepsilon^2\delta^2),\\
& & \frac{\partial \Delta p_{nh} }{\partial z} = - \varepsilon^2\delta \frac{\partial w}{\partial t} + {\cal O}(\varepsilon^2\delta^2),
\end{eqnarray*}
leading to
\begin{eqnarray*}
{\cal C}_1 & = & \frac{\partial}{\partial x} \int_{z_b}^{\delta\eta} u \Delta p_{nh} + \varepsilon^2\delta\int_{z_b}^{\delta\eta} w \frac{\partial w}{\partial t} - \frac{\partial z_b}{\partial t} \left. \Delta p_{nh}\right|_b + {\cal O}(\varepsilon^2\delta^2) \nonumber\\
& = & \frac{\partial}{\partial x} \int_{z_b}^{\delta\eta} u \Delta p_{nh} + \varepsilon^2\delta\frac{\partial }{\partial t}\int_{z_b}^{\delta\eta} \frac{w^2}{2} - \frac{\partial z_b}{\partial t} \left. \Delta p_{nh}\right|_b + {\cal O}(\varepsilon^2\delta^2). \nonumber
\end{eqnarray*}
Due to the rescaling applied to the time derivative of $z_b$ (see
paragraph~\ref{subsec:NS_rescaled}), the Leibniz rule applied to obtain
the preceding relation reads
$$\int_{z_b}^{\delta\eta} w \frac{\partial
w}{\partial t} =  \delta\frac{\partial z_b}{\partial
t}\frac{w_b^2}{2} - \delta\frac{\partial \eta}{\partial
t}\frac{w_s^2}{2} + \frac{\partial }{\partial
t}\int_{z_b}^{\delta\eta} \frac{w^2}{2} = \frac{\partial }{\partial
t}\int_{z_b}^{\delta\eta} \frac{w^2}{2}+ {\cal O}(\delta).$$
And finally we have for ${\cal C}_1$
\begin{eqnarray*}
{\cal C}_1 & = & \frac{\partial}{\partial x} \int_{z_b}^{\delta\eta} u \Delta p_{nh} + \varepsilon^2\delta\frac{\partial }{\partial t}\int_{z_b}^{\delta\eta} \frac{w^2}{2} - \frac{\partial z_b}{\partial t} \left. \Delta p_{nh}\right|_b + {\cal O}(\varepsilon^2\delta^2). \nonumber
\end{eqnarray*}
From relations (\ref{eq:u_approx}) and (\ref{eq:u_bar}) we obtain
$$u = \left(1 +
\frac{\varepsilon\kappa_0}{\nu_0}\left(z-z_b-\frac{(z-z_b)^2}{2H_\delta}-\frac{H_\delta}{3}\right)\right)\bar{u}
= \left(1 + \varepsilon f(z-z_b,H_\delta)\right)\bar{u} + {\cal O}(\varepsilon^2),$$
so we have for ${\cal C}_2$ and ${\cal C}_3$
\begin{eqnarray*}
{\cal C}_2 & = & \frac{\partial}{\partial x}\int_{z_b}^{\delta\eta} \bar{u}\left(2\varepsilon\delta\nu_0\frac{\partial u}{\partial x}\right) - 2\varepsilon\delta\nu_0\int_{z_b}^{\delta\eta} \frac{\partial \bar{u}}{\partial x}\left(\frac{\partial u}{\partial x}\right),\\
& = & \frac{\partial}{\partial x}\int_{z_b}^{\delta\eta} \bar{u}\left(2\varepsilon\delta\nu_0\frac{\partial u}{\partial x}\right) - 2\varepsilon\delta\nu_0\left(\int_{z_b}^{\delta\eta} \left(\frac{\partial \bar{u}}{\partial x}\right)^2 + \varepsilon\frac{\partial \bar{u}}{\partial x}\int_{z_b}^{\delta\eta} \frac{\partial (f\bar{u})}{\partial x}\right) + {\cal O}(\varepsilon^3\delta),\\
{\cal C}_3 & = & \frac{\partial}{\partial x}\int_{z_b}^{\delta\eta} \bar{u} p_{v,f} - \int_{z_b}^{\delta\eta} \frac{\partial \bar{u}}{\partial x} p_{v,f},\\
& = & \frac{\partial}{\partial x}\int_{z_b}^{\delta\eta} \bar{u} p_{v,f} + \int_{z_b}^{\delta\eta} \frac{\partial w}{\partial z} p_{v,f} + \varepsilon\int_{z_b}^{\delta\eta} \frac{\partial (f\bar{u})}{\partial x} p_{v,f} + {\cal O}(\varepsilon^3\delta),\\
& = & \frac{\partial}{\partial x}\int_{z_b}^{\delta\eta} \bar{u} p_v + [w p_{v,f}]_{z_b}^{\delta\eta} - \varepsilon\delta\nu_0\int_{z_b}^{\delta\eta} w\frac{\partial^2 u}{\partial x\partial z} - 2\varepsilon\delta\nu_0\int_{z_b}^{\delta\eta} w\frac{\partial^2 w}{\partial z^2}\\
& & + \varepsilon\int_{z_b}^{\delta\eta} \frac{\partial (f\bar{u})}{\partial x} p_{v,f} + {\cal O}(\varepsilon^3\delta),\\
& = & \frac{\partial}{\partial x}\int_{z_b}^{\delta\eta} \bar{u} p_{v,f} + [w p_{v,f}]_{z_b}^{\delta\eta} - 2\varepsilon\delta\nu_0\int_{z_b}^{\delta\eta} \frac{\partial}{\partial z}\left(w\frac{\partial w}{\partial z}\right) + 2\varepsilon\delta\nu_0\int_{z_b}^{\delta\eta} \left(\frac{\partial w}{\partial z}\right)^2\\
& & - \varepsilon\delta\nu_0\frac{\partial}{\partial x}\int_{z_b}^{\delta\eta} \left(w\frac{\partial u}{\partial z}\right) + \varepsilon\delta^2\nu_0\frac{\partial \eta}{\partial x}w_s\left.\frac{\partial u}{\partial z}\right|_s - \varepsilon\delta\nu_0\frac{\partial z_b}{\partial x}w_b\left.\frac{\partial u}{\partial z}\right|_b \nonumber\\
& & + \varepsilon\delta\nu_0\int_{z_b}^{\delta\eta} \frac{\partial w}{\partial x}\frac{\partial u}{\partial z} + \varepsilon\int_{z_b}^{\delta\eta} \frac{\partial (f\bar{u})}{\partial x} p_{v,f} + {\cal O}(\varepsilon^3\delta),\\
& = & \frac{\partial}{\partial x}\int_{z_b}^{\delta\eta} \bar{u} p_{v,f} + [w p_{v,f}]_{z_b}^{\delta\eta} - 2\varepsilon\delta\nu_0\left( w_s \left.\frac{\partial w}{\partial z}\right|_s - w_b \left.\frac{\partial w}{\partial z}\right|_b \right) \\
& & + 2\varepsilon\delta\nu_0\int_{z_b}^{\delta\eta} \left(\frac{\partial w}{\partial z}\right)^2 - \varepsilon\delta\nu_0\frac{\partial}{\partial x}\int_{z_b}^{\delta\eta} \left(w\frac{\partial u}{\partial z}\right) + \varepsilon\delta^2\nu_0\frac{\partial \eta}{\partial x}w_s\left.\frac{\partial u}{\partial z}\right|_s \\
& & - \varepsilon\delta\nu_0\frac{\partial z_b}{\partial x}w_b\left.\frac{\partial u}{\partial z}\right|_b + \varepsilon\delta\nu_0\int_{z_b}^{\delta\eta} \frac{\partial w}{\partial x}\frac{\partial u}{\partial z} + \varepsilon\int_{z_b}^{\delta\eta} \frac{\partial (f\bar{u})}{\partial x} p_{v,f} + {\cal O}(\varepsilon^3\delta),
\end{eqnarray*}
and from relation (\ref{eq:duzb}) we also have
\begin{eqnarray*}
\nu_0\int_{z_b}^{\delta\eta} \frac{\partial w}{\partial x}\frac{\partial u}{\partial z} & = & \varepsilon\kappa_0\int_{z_b}^{\delta\eta} \frac{\partial w}{\partial x}\left(1 - \frac{z - z_b}{H_\delta}\right)u_b + {\cal O}(\varepsilon^2),\\
& = & \varepsilon\kappa_0 \frac{H_\delta}{2}\frac{\partial^2 z_b}{\partial x\partial t}u_b + \varepsilon\kappa_0\int_{z_b}^{\delta\eta} \left(-z\frac{\partial^2 \bar{u}}{\partial x^2} + \frac{\partial^2 (z_b\bar{u})}{\partial x^2} \right)\left(1 - \frac{z - z_b}{H_\delta}\right)\bar{u} + {\cal O}(\varepsilon^2),\\
& = & \varepsilon\kappa_0 \frac{H_\delta}{2}\frac{\partial^2 z_b}{\partial x\partial t}u_b - \varepsilon\kappa_0\frac{H_\delta^2}{6}\frac{\partial^2 \bar{u}}{\partial x^2}\bar{u} + \varepsilon\kappa_0\frac{H_\delta}{2}\frac{\partial}{\partial x}\left(\frac{\partial z_b}{\partial x}\bar{u}^2\right) + {\cal O}(\varepsilon^2),\\
& = & \varepsilon\kappa_0 \frac{H_\delta}{2}\frac{\partial^2 z_b}{\partial x\partial t}u_b - \varepsilon\kappa_0\frac{\partial}{\partial x}\left(\frac{H_\delta^2}{6}\frac{\partial \bar{u}}{\partial x}\bar{u}\right) + \varepsilon\kappa_0\frac{H_\delta^2}{6}\left(\frac{\partial \bar{u}}{\partial x}\right)^2\\
& & + \varepsilon\kappa_0 \frac{H_\delta}{3}\frac{\partial H_\delta}{\partial x}\frac{\partial \bar{u}}{\partial x}\bar{u} + \varepsilon\kappa_0\frac{\partial}{\partial x}\left(\frac{H_\delta}{2}\frac{\partial z_b}{\partial x}\bar{u}^2\right) - \frac{\varepsilon\kappa_0}{2}\frac{\partial H_\delta}{\partial x}\frac{\partial z_b}{\partial x}\bar{u}^2 + {\cal O}(\varepsilon^2).
\end{eqnarray*}
The preceding expression shows that due to relation (\ref{eq:duzb}), the
term
$$\nu_0\int_{z_b}^{\delta\eta} \frac{\partial w}{\partial
x}\frac{\partial u}{\partial z},$$
has to be treated as a friction term in the energy equality. We finally have for ${\cal R} = {\cal C}_2 - {\cal C}_3 + {\cal C}_4$
\begin{eqnarray}
{\cal R} & = & \frac{\partial}{\partial x}\int_{z_b}^{\delta\eta} \left(2\varepsilon\delta\nu_0\bar{u}\frac{\partial u}{\partial x} -p_{v,f}\right) - 2\varepsilon\delta\nu_0\int_{z_b}^{\delta\eta} \left(\left(\frac{\partial \bar{u}}{\partial x}\right)^2 + \left(\frac{\partial w}{\partial z}\right)^2\right) \nonumber\\
& & + \varepsilon\delta\nu_0\frac{\partial}{\partial x}\int_{z_b}^{\delta\eta} \left(w\frac{\partial u}{\partial z}\right) + \varepsilon^2\delta\frac{\partial}{\partial x}\left(\kappa_0\frac{H_\delta}{2}\left(\frac{H_\delta}{3}\frac{\partial \bar{u}}{\partial x}\bar{u} - \frac{\partial z_b}{\partial x}\bar{u}^2\right)\right)\nonumber\\
& & - \frac{\varepsilon^2\delta\kappa_0}{6}\left(\left( H_\delta\frac{\partial \bar{u}}{\partial x} + \frac{\partial H_\delta}{\partial x}\bar{u}\right)^2 - \left(-2\left(\frac{\partial z_b}{\partial x}\right)^2 + \delta\frac{\partial \eta}{\partial x}\frac{\partial z_b}{\partial x}\right.\right.\nonumber\\
& & + \left.\left. \delta^2\left(\frac{\partial \eta}{\partial x}\right)^2\right)\bar{u}^2\right) + \frac{\partial z_b}{\partial t}\left(\left. p_{v,f}\right|_b + 2\varepsilon\delta\nu_0\left.\frac{\partial u}{\partial x}\right|_b \right)\nonumber\\
& & - \varepsilon^2\delta\kappa_0 \frac{H_\delta}{2}\frac{\partial^2 z_b}{\partial x\partial t}u_b + {\cal O}(\varepsilon^3\delta). \nonumber
\end{eqnarray}
Returning to the initial variables and integrating ${\cal C}_1$ and
${\cal R}$ into relation (\ref{eq:energy}) gives an energy equality for the model (\ref{eq:sv_nonhyd_init1})-(\ref{eq:sv_nonhyd_init2}) under the form
\begin{eqnarray*}
\lefteqn{\frac{\partial}{\partial t} \left( E_h + \frac{H\overline{w^2}}{2}\right) + \frac{\partial}{\partial x}\left(\bar{u}\left( E_h + H \bar{p}_{nh}\right) - \nu\int_{z_b}^{\eta} \left(2H\bar{u}\frac{\partial u}{\partial x} + w\frac{\partial u}{\partial z}\right) \right)} \\
& & \quad + \frac{\partial}{\partial x}\left(\kappa\left(\frac{z_b^2}{6}\frac{\partial \bar{u}}{\partial x}\bar{u} + \frac{z_b}{2}\frac{\partial z_b}{\partial x}\bar{u}^2\right)\right) \nonumber\\
& & = - 2\nu\int_{z_b}^{\eta} \left(\left(\frac{\partial \bar{u}}{\partial x}\right)^2 + \left(\frac{\partial u}{\partial x}\right)^2\right)  - \frac{\kappa}{6}\left( z_b\frac{\partial \bar{u}}{\partial x} + \frac{\partial z_b}{\partial x}\bar{u}\right)^2\nonumber\\
& & \quad - \frac{\kappa}{1+\frac{\kappa H}{3 \nu}}\left(1+\frac{11}{6}\left(\frac{\partial z_b}{\partial x}\right)^2\right) \bar{u}^2\nonumber\\
& & \quad - H\frac{\partial p^a}{\partial t} + \left(\left. p_{nh}\right|_b + 2\nu\frac{\partial u_b}{\partial x}\right)\frac{\partial z_b}{\partial t} + \kappa\frac{z_b}{2}\frac{\partial^2 z_b}{\partial x\partial t}\bar{u},\nonumber
\end{eqnarray*}
where
\begin{eqnarray*}
H\overline{w^2}& = & \int_{z_b}^{\eta} w^2 = \int_{z_b}^{\eta} \left(\frac{\partial z_b}{\partial t} - z\frac{\partial \bar{u}}{\partial x} + \frac{\partial (z_b\bar{u})}{\partial x} \right)^2 = -\frac{z_b^2}{3}\left(\frac{\partial \bar{u}}{\partial x}\right)^2\nonumber\\
& & - z_b\left(\frac{\partial (z_b\bar{u})}{\partial x}\right)^2 + z_b^2\frac{\partial \bar{u}}{\partial x}\frac{\partial (z_b\bar{u})}{\partial x} - z_b\left(\frac{\partial z_b}{\partial t}\right)^2\nonumber\\
& & + 2z_b\frac{\partial z_b}{\partial t}\left(\frac{z_b}{2}\frac{\partial \bar{u}}{\partial x} - \frac{\partial (z_b\bar{u})}{\partial x}\right),\\
H \bar{p}_{nh} & = & \int_{z_b}^{\eta} p_{nh}.
\end{eqnarray*}
When the time derivatives of $p^a$ and $z_b$ are dropped, the right hand
side of the preceding energy equality is always negative.

\subsection{A more complex approximation, $\delta = {\cal O}(1)$}
\label{subsec:SV_nhyd2}

Now we return to the dimensionless and rescaled variables. The
assumption that the elevation of the free surface is small done in
paragraph~\ref{subsec:nhyd1} is now relaxed i.e. $\delta = {\cal O}(1)$. This means that no assumption is made concerning the hydraulic regime. We consider for the
pressure $p$ the complete expression obtained in
(\ref{eq:pfin}) and the following proposition is a refinement of the Proposition~\ref{prop:SV+}.
\begin{proposition}
The system defined by
\begin{eqnarray}
\lefteqn{\frac{\partial H}{\partial t} + \frac{\partial}{\partial x} \bigl(H\bar{u}\bigr) = 0,} \label{eq:sv++_1}\\
\lefteqn{\frac{\partial}{\partial t}(H \bar{u}) + \frac{\partial}{\partial x} (H_m \bar{u}^2) + \frac{1}{2}\frac{\partial H^2}{\partial x} + \frac{\partial (H\bar{p}_{ng,nv})}{\partial x} = - H\frac{\partial p^a}{\partial x} - gH\frac{\partial z_b}{\partial x} \nonumber}\\
& & + \frac{\partial }{\partial x}\bigl(4\nu H\frac{\partial \bar{u}}{\partial x}\bigr) + \frac{\partial}{\partial x}\left(\kappa({\bf v}_b,H) H\left( \frac{H}{6}\frac{\partial \bar{u}}{\partial x} - \left(\frac{7}{6}\frac{\partial z_b}{\partial x} + \frac{1}{3}\frac{\partial \eta}{\partial x}\right)\bar{u} \right)\right) \nonumber\\
& & -\frac{\kappa({\bf v}_b,H)}{1 + \frac{\kappa({\bf v}_b,H) H}{3 \nu}}\left(1 + \frac{5}{2}\left(\frac{\partial z_b}{\partial x}\right)^2\right)\bar{u} + \kappa({\bf v}_b,H)\frac{\partial z_b}{\partial x}\left(\bigl(\frac{1}{2}\frac{\partial H}{\partial x}+\frac{\partial z_b}{\partial x}\bigr)\bar{u}\right. \nonumber\\
& & \left. + \frac{H}{2}\frac{\partial \bar{u}}{\partial x}\right) - \frac{\partial z_b}{\partial x}\left. p_{ng,nv}\right|_b  + z_b\frac{\partial z_b}{\partial x}\frac{\partial^2 z_b}{\partial t^2} - \frac{1}{2}\frac{\partial}{\partial x}\left(H^2\frac{\partial^2 z_b}{\partial t^2}\right),\label{eq:sv++_2}
\end{eqnarray}
where $\bar{{\bf v}}=(1,\frac{\partial z_b}{\partial x})^T\frac{\bar{u}}{1+\frac{\kappa_l H}{3\nu}}$
results from an approximation in ${\cal O}(\varepsilon^3\delta)$ of the Navier-Stokes equations.
In the previous expressions, $H_m$ is a modified water height taking
 into account the Coriolis-Boussinesq coefficient and $\bar{p}_{ng,nv}$,
 $\left. p_{ng,nv}\right|_b$ corresponds to the vertically averaged and bottom
 value of the non gravitational and non viscous part of the pressure $p$ given by (\ref{eq:pfin}).\label{prop:SV++}
\end{proposition}

{\it Proof of prop.~\ref{prop:SV++}:} we still start from the averaged momentum equation~(\ref{eq:momentum})
where, compared to the first extension of the Saint-Venant model
detailed in paragraphs~\ref{subsec:extension_1}, \ref{subsec:nhyd1} and
\ref{subsec:prop}, the expressions of
$$\int_{z_b}^{\delta\eta}p\ dz\quad\mbox{and}\quad \int_{z_b}^{\delta\eta}u^2 dz,$$
have to be refined. The approximation $\overline{u^2} = \bar{u}^2 +
{\cal O}(\varepsilon^2)$ obtained in paragraph~\ref{subsec:SV} is no
more sufficient. From (\ref{eq:cons_w}), (\ref{eq:phyd}),
(\ref{eq:u_approx}) and (\ref{eq:u_bar}) we get
\begin{eqnarray}
u & = & \left(1 +
\frac{\varepsilon\kappa_0}{\nu_0}\left(z-z_b-\frac{(z-z_b)^2}{2H_\delta}-\frac{H_\delta}{3}\right)\right)\bar{u} + {\cal O}(\varepsilon^2) \nonumber\\
& = & \bigl(1+\varepsilon f(z-z_b,H_\delta)\bigr)\bar{u} + {\cal O}(\varepsilon^2),\nonumber\\
w & = & \frac{\partial z_b}{\partial t} - \frac{\partial}{\partial x}\left(\left(z-z_b+\varepsilon\int_{z_b}^z f(z,z_b,H_\delta) dz\right)\bar{u}\right) + {\cal O}(\varepsilon^2)\nonumber\\
& = & \frac{\partial z_b}{\partial t} - \frac{\partial}{\partial x}\bigl(g(z,z_b,H_\delta)\bar{u}\bigr) + {\cal O}(\varepsilon^2),\nonumber\\
\frac{\partial p}{\partial x} & = & \delta\frac{\partial}{\partial x}(p^a + \eta) - 2\varepsilon\delta\frac{\partial}{\partial x}\left(\nu_0 \frac{\partial \bar{u}}{\partial x}\right) + {\cal O}(\varepsilon^2\delta),\nonumber
\end{eqnarray}
and Eq.~(\ref{eq:final1_2}) is equivalent to
\begin{eqnarray*}
\delta\frac{\partial \bar{u}}{\partial t} + \delta^2\bar{u}\frac{\partial \bar{u}}{\partial x} + \delta\frac{\partial \eta}{\partial x} & = & -\delta\frac{\partial }{\partial x}p^a - \frac{\delta\kappa_0}{H_\delta\bigl(1 + \frac{\varepsilon\kappa_0}{3\nu_0}H_\delta\bigr)} \bar{u} \nonumber\\
& & + \frac{\varepsilon\delta}{H_\delta}\frac{\partial}{\partial x}\bigl(4 \nu_0 H_\delta\frac{\partial \bar{u}}{\partial x}\bigr) + {\cal O}(\varepsilon^2\delta).
\end{eqnarray*}
Now we can improve the approximation (\ref{eq:u_est}) in the following way
\begin{eqnarray*}
\delta\frac{\partial }{\partial z}\bigl(\nu_0\frac{\partial u}{\partial z}\bigr) & = & \varepsilon\delta\frac{\partial u}{\partial t} + \varepsilon\delta^2 u\frac{\partial u}{\partial x} + \varepsilon\delta^2 w\frac{\partial u}{\partial z} + \varepsilon\frac{\partial p}{\partial x} - \varepsilon^2\delta\frac{\partial }{\partial x}\left(\nu_0\frac{\partial u}{\partial x}\right)\nonumber\\
& = & \varepsilon\delta\frac{\partial}{\partial t}\bigl((1+\varepsilon f)\bar{u}\bigr) + \varepsilon\delta^2 (1+\varepsilon f)\bar{u}\frac{\partial}{\partial x}\bigl((1+\varepsilon f)\bar{u}\bigr) \nonumber\\
& & + \varepsilon\delta\frac{\partial }{\partial x}(\eta + p^a) + \varepsilon\delta^2 w\frac{\partial u}{\partial z} - 3\varepsilon^2\delta\frac{\partial}{\partial x}\left(\nu_0 \frac{\partial \bar{u}}{\partial x}\right) + {\cal O}(\varepsilon^3\delta)\nonumber\\
& = & - \frac{\varepsilon\delta\kappa_0}{H_\delta\bigl(1 + \frac{\varepsilon\kappa_0}{3\nu_0}H_\delta\bigr)} \bar{u} + \varepsilon\delta^2 w\frac{\partial u}{\partial z}  + \varepsilon^2\delta\frac{\partial}{\partial t}\bigl(f\bar{u}\bigr) + \varepsilon^2\delta^2 \frac{\partial}{\partial x}\bigl(f\bar{u}^2\bigr)\nonumber\\
& & + \varepsilon^2\delta\frac{\partial}{\partial x}\left(\nu_0 \frac{\partial \bar{u}}{\partial x}\right) + \varepsilon^2\delta\frac{4\nu_0}{H_\delta}\frac{\partial H_\delta}{\partial x}\frac{\partial \bar{u}}{\partial x} + {\cal O}(\varepsilon^3\delta).\nonumber
\end{eqnarray*}
Taking into acount the boundary condition (\ref{eq:BC_ap5}), an integration of the preceding relation from $z_b$ to $z$ gives
\begin{eqnarray}
\nu_0\frac{\partial u}{\partial z} & = & \frac{\varepsilon\kappa_0}{\bigl(1 + \frac{\varepsilon\kappa_0}{3\nu_0}H_\delta}\bigr)\left(1 - \frac{z-z_b}{H_\delta}\right) \bar{u} + \varepsilon^2\delta\left(\delta\bar{u}fw + \delta\bar{u}\frac{\partial\bar{u}}{\partial x}\int_{z_b}^z f \right.\nonumber\\
& & \left. + \frac{\partial}{\partial t}\left(\bar{u}\int_{z_b}^z f\right)+ \delta\frac{\partial}{\partial x}\left(\bar{u}^2\int_{z_b}^z f \right) + (z-z_b)\frac{\partial}{\partial x}\left(\nu_0 \frac{\partial \bar{u}}{\partial x}\right)\right. \nonumber\\
& & \left. + \frac{4\nu_0(z-z_b)}{H_\delta}\frac{\partial H_\delta}{\partial x}\frac{\partial \bar{u}}{\partial x}\right) + {\cal O}(\varepsilon^3\delta),\label{eq:u22}
\end{eqnarray}
where the relation
$$\int_{z_b}^z w \frac{\partial u}{\partial z} =
\varepsilon\bar{u}\int_{z_b}^z w \frac{\partial f}{\partial z} =
\varepsilon\bar{u} (fw -f|_b w|_b) + \varepsilon\bar{u}\frac{\partial
\bar{u}}{\partial x}\int_{z_b}^z f,$$
has been used. Another integration of relation (\ref{eq:u22}) between
$z_b$ and $z$ gives
\begin{eqnarray}
u & = & \left(1+\frac{\varepsilon\kappa_0}{\nu_0}\left(z-z_b - \frac{(z-z_b)^2}{2H_\delta}\right)\right)u_b + \frac{\varepsilon^2\delta^2}{\nu_0}\bar{u}\int_{z_b}^z fw\nonumber\\
& & + \frac{\varepsilon^2\delta^2}{\nu_0}\bar{u}\frac{\partial\bar{u}}{\partial x}\int_{z_b}^z\int_{z_b}^{z_1} f + \frac{\varepsilon^2\delta}{\nu_0}\frac{\partial}{\partial t}\left(\bar{u}\int_{z_b}^z\int_{z_b}^{z_1} f\right) + \frac{\varepsilon^2\delta^2}{\nu_0} \frac{\partial}{\partial x}\left(\bar{u}^2\int_{z_b}^z\int_{z_b}^{z_1} f \right)\nonumber\\
& & + \frac{\varepsilon^2\delta}{2} (z-z_b)^2\frac{\partial^2 \bar{u}}{\partial x^2} + \varepsilon^2\delta\frac{2(z-z_b)^2}{H_\delta}\frac{\partial H_\delta}{\partial x}\frac{\partial \bar{u}}{\partial x} + {\cal O}(\varepsilon^3\delta),\nonumber\\
& = & \left(1+\frac{\varepsilon\kappa_0}{\nu_0}\left(z-z_b - \frac{(z-z_b)^2}{2H_\delta}\right)\right)u_b + \varepsilon^2\delta \Delta u+ {\cal O}(\varepsilon^3\delta),\nonumber
\end{eqnarray}
so we obtain the new expressions for $\bar{u}$, $\bar{u}^2$ and $u^2$
\begin{eqnarray}
\bar{u} & = & \left(1+\frac{\varepsilon\kappa_0}{3\nu_0}H_\delta\right)u_b + \varepsilon^2\delta \overline{\Delta u}+ {\cal O}(\varepsilon^3\delta),\nonumber\\
\bar{u}^2 & = & \left(1+\frac{2\varepsilon\kappa_0}{3\nu_0}H_\delta\right)u_b^2 + 2\varepsilon^2\delta \overline{\Delta u} + {\cal O}(\varepsilon^3\delta),\nonumber\\
u^2 & = & \left(1+\frac{2\varepsilon\kappa_0}{\nu_0}\left(z-z_b - \frac{(z-z_b)^2}{2H_\delta}\right) +\frac{\varepsilon^2\kappa_0^2}{\nu_0}\left(z-z_b - \frac{(z-z_b)^2}{2H_\delta}\right)^2\right)u_b\nonumber\\
& & + 2\varepsilon^2\delta \Delta u + {\cal O}(\varepsilon^3\delta),\nonumber
\end{eqnarray}
so finally
\begin{eqnarray*}
\overline{u^2} & = & \left(1+\frac{2\varepsilon\kappa_0}{3\nu_0}H_\delta + \frac{2\varepsilon^2\kappa_0^2}{15\nu_0^2}H_\delta^2 \right)u_b^2 + 2\varepsilon^2\delta \overline{\Delta u}\nonumber\\
 & = & \left(1 + \frac{2\varepsilon^2\kappa_0^2}{15\nu_0^2}H_\delta^2 \right)\bar{u}^2 + {\cal O}(\varepsilon^3\delta).\nonumber
\end{eqnarray*}
Now concerning the expression of the pressure terms, it has to be noticed that Eqs.~(\ref{eq:pfin}) and (\ref{eq:pres}) only
differ by the terms
$${\cal A} = \varepsilon^2\delta^2\frac{\partial }{\partial x}\int_z^{\delta
\eta} uw\ dz - \varepsilon^2\delta^2 w^2.$$
Using
$$u=\bar{u}+ {\cal O}(\varepsilon), \quad w = \frac{\partial
z_b}{\partial t} - \frac{\partial}{\partial x} \int_{z_b}^z u\ dz,$$
it comes
\begin{eqnarray*}
{\cal A} & = & \varepsilon^2\delta^2\left (-\frac{\delta^2\eta^2 - z^2}{2} \frac {\partial }{\partial x}\bigl(\bar{u}\frac {\partial \bar{u}}{\partial x}\bigr) - \delta^2\eta\frac {\partial \eta}{\partial x}\frac {\partial \bar{u}}{\partial x}\bar{u} + (\delta \eta -z)\frac{\partial}{ \partial x}\bigl(\bar{u}\frac {\partial (z_b \bar{u})}{\partial x}\bigr)\right.\\
& & + \left.\delta \frac{\partial \eta}{\partial x} \frac{\partial (z_b \bar{u})}{\partial x}\bar{u} -\left(- z\frac{\partial \bar{u}}{\partial x} + \frac{\partial (z_b\bar{u})}{\partial x}\right)^2\right).
\end{eqnarray*}
This leads to the new expression for the fluid pressure $p$ appearing in
(\ref{eq:momentum})
\begin{eqnarray*}
\int_{z_b}^{\delta \eta} p\ dz & = & \int_{z_b}^{\delta \eta} (p_{nh} + {\cal A}) dz\\
& = & \int_{z_b}^{\delta \eta} p_{nh} dz + \frac{\varepsilon^2\delta^2 H_\delta}{6}\left( -4H_\delta^2\bigl(\frac{\partial \bar{u}}{\partial x}\bigr)^2 - 2 H_\delta^2\bar{u}\frac{\partial^2 \bar{u}}{\partial x^2} -6 H_\delta\frac{\partial H_\delta}{\partial x}\frac{\partial \bar{u}}{\partial x}\bar{u}\right. \nonumber\\
& & \left. +9 H_\delta\frac{\partial z_b}{\partial x}\frac{\partial \bar{u}}{\partial x}\bar{u} + 3H_\delta\frac{\partial^2 z_b}{\partial x^2}\bar{u}^2 + 6\frac{\partial z_b}{\partial x}\frac{\partial H_\delta}{\partial x}\bar{u}^2\right) + {\cal O}(\varepsilon^3\delta),
\end{eqnarray*}
where $\int_{z_b}^{\delta \eta} p_{nh}=\int_{z_b}^{\delta \eta} (p_h +
\Delta p_{nh})$ is given by (\ref{eq:presint}). Conversely using (\ref{eq:pfin}) we obtain
\begin{eqnarray*}
p_b = \left. p_{nh}\right|_b + \frac{\varepsilon^2\delta^2}{2}\left(-\frac{\partial}{\partial x}\bigl(H_\delta^2 \frac{\partial \bar{u}}{\partial x}\bar{u}\bigr) + 4 H_\delta\frac{\partial z_b}{\partial x}\frac{\partial \bar{u}}{\partial x}\bar{u} + 2\frac{\partial}{\partial x}\bigl(H_\delta\frac{\partial z_b}{\partial x}\bigr)\bar{u}^2\right) + {\cal O}(\varepsilon^3\delta),
\end{eqnarray*}
where $\left. p_{nh}\right|_b$ is given by (\ref{eq:presint2}). Inserting (\ref{eq:du}), (\ref{eq:pb}) and (\ref{eq:pm}) in equilibrium
(\ref{eq:momentum}) leads to the system
\begin{eqnarray}
\lefteqn{\frac{\partial H_\delta}{\partial t} + \frac{\partial}{\partial x} \bigl(H_\delta\bar{u}\bigr) = 0,}\label{eq:sv_nonhyd_fin1}\\
\lefteqn{\delta\frac{\partial}{\partial t}(H_\delta \bar{u}) + \delta^2\frac{\partial}{\partial x} (H_{\delta,m} \bar{u}^2) + \frac{1}{2}\frac{\partial H_\delta^2}{\partial x} + \frac{\partial (H_\delta\bar{p}_{ng,nv})}{\partial x} = - H_\delta\frac{\partial }{\partial x}(\delta p^a + z_b) }\nonumber\\
& & + \varepsilon\delta\frac{\partial }{\partial x}\bigl(4\nu_0H_\delta \frac{\partial \bar{u}}{\partial x}\bigr) + \varepsilon^2\delta\frac{\partial}{\partial x}\left(\kappa_0 H_\delta\bigl( \frac{H_\delta}{6}\frac{\partial \bar{u}}{\partial x} - \bigl(\frac{7}{6}\frac{\partial z_b}{\partial x} + \frac{\delta}{3}\frac{\partial \eta}{\partial x}\bigr)\bar{u} \bigr)\right) \nonumber\\
& & + \varepsilon^2\delta\kappa_0\frac{\partial z_b}{\partial x}\left(\bigl(\frac{1}{2}\frac{\partial H_\delta}{\partial x}+\frac{\partial z_b}{\partial x}\bigr)\bar{u} + \frac{H_\delta}{2}\frac{\partial \bar{u}}{\partial x}\right) - \frac{\partial z_b}{\partial x} \left. p_{ng,nv}\right|_b\nonumber\\
& & -\delta\kappa_0\left(1 + \frac{5\varepsilon^2}{2}\left(\frac{\partial z_b}{\partial x}\right)^2\right)u_b +\varepsilon^2\delta z_b\frac{\partial z_b}{\partial x}\frac{\partial^2 z_b}{\partial t^2} - \frac{1}{2}\frac{\partial}{\partial x}\left(H_\delta^2\frac{\partial^2 z_b}{\partial t^2}\right) + {\cal O}(\varepsilon^3\delta), \label{eq:sv_nonhyd_fin2}
\end{eqnarray}
where 
\begin{eqnarray*}
H_\delta^m & = & H_\delta  \left(1 + \frac{2\varepsilon^2\kappa_0^2}{15\nu_0^2}H_\delta^2 \right),\nonumber\\ 
H_\delta \bar{p}_{ng,nv} & = & \int_{z_b}^{\delta\eta} (p - p_h)\ dz\nonumber\\
& = & \varepsilon^2\delta\frac{\partial}{\partial x}\left(\frac{H_\delta^3}{6}\frac{\partial^2 \bar{u}}{\partial x \partial t} + \frac{H_\delta^2}{2}\frac{\partial^2 (z_b\bar{u})}{\partial x \partial t} - \delta\eta\frac{H_\delta^2}{2}\frac{\partial^2 \bar{u}}{\partial x \partial t}\right.\\
& & \left.  - \delta H_\delta \frac{\partial \eta}{\partial t} \bigl(\delta\eta\frac{\partial \bar{u}}{\partial x}- \frac{\partial (z_b\bar{u})}{\partial x} \bigr)\right) + \frac{\varepsilon^2\delta^2 H_\delta}{6}\left( -4H_\delta^2\bigl(\frac{\partial \bar{u}}{\partial x}\bigr)^2\right.\\
& & \left. - 2 H_\delta^2\bar{u}\frac{\partial^2 \bar{u}}{\partial x^2} -6 H_\delta\frac{\partial H_\delta}{\partial x}\frac{\partial \bar{u}}{\partial x}\bar{u} + 9 H_\delta\frac{\partial z_b}{\partial x}\frac{\partial \bar{u}}{\partial x}\bar{u} \right. \nonumber\\
 & & \left. + 3H_\delta\frac{\partial^2 z_b}{\partial x^2}\bar{u}^2 + 6\frac{\partial z_b}{\partial x}\frac{\partial H_\delta}{\partial x}\bar{u}^2\right) + \varepsilon^2\delta\frac{H_\delta^2}{2}\frac{\partial^2 z_b}{\partial t^2} + {\cal O}(\varepsilon^3\delta),
\end{eqnarray*}
and
\begin{eqnarray*}
\left.p_{ng,nv}\right|_b & = & \left. (p - p_h)\right|_b\\
& = & \frac{\varepsilon^2\delta}{2}\left( -\frac{\partial}{\partial t}\bigl(H_\delta^2\frac{\partial \bar{u}}{\partial x}\bigr) + 2H_\delta\frac{\partial}{\partial t}\bigl(\frac{\partial z_b}{\partial x}\bar{u}\bigr) +2\delta\frac{\partial \eta}{\partial t}\frac{\partial z_b}{\partial x}\bar{u}\right)\nonumber\\
& & + \varepsilon^2\delta\left(H_\delta\frac{\partial^2 z_b}{\partial t^2} + \delta\frac{\partial \eta}{\partial t}\frac{\partial z_b}{\partial t}\right) + 4 H_\delta\frac{\partial z_b}{\partial x}\frac{\partial \bar{u}}{\partial x}\bar{u} + 2\frac{\partial}{\partial x}\bigl(H_\delta\frac{\partial z_b}{\partial x}\bigr)\bar{u}^2\nonumber\\
& & + \frac{\varepsilon^2\delta^2}{2}\left(-\frac{\partial}{\partial x}\bigl(H_\delta^2 \frac{\partial \bar{u}}{\partial x}\bar{u}\bigr)\right)  + \varepsilon^2\delta H_\delta\frac{\partial^2 z_b}{\partial t^2} + {\cal O}(\varepsilon^3\delta),\nonumber\\
& = & \frac{\varepsilon^2\delta}{2}\left( H_\delta^2\frac{\partial^2 \bar{u}}{\partial x\partial t} + 2H_\delta\frac{\partial^2 (z_b\bar{u})}{\partial x\partial t} +2\delta\frac{\partial \eta}{\partial t}\frac{\partial (z_b\bar{u})}{\partial x} - 2\delta\eta\bigl(\delta\frac{\partial \eta}{\partial t}\frac{\partial \bar{u}}{\partial x}\right.\\
& & \left. + H_\delta\frac{\partial^2 \bar{u}}{\partial x\partial t}\bigr)\right) + \varepsilon^2\delta\left(H_\delta\frac{\partial^2 z_b}{\partial t^2} + \delta\frac{\partial \eta}{\partial t}\frac{\partial z_b}{\partial t}\right) + \varepsilon^2\delta H_\delta\frac{\partial^2 z_b}{\partial t^2}\nonumber\\
& & + \frac{\varepsilon^2\delta^2}{2}\left(-\frac{\partial}{\partial x}\bigl(H_\delta^2 \frac{\partial \bar{u}}{\partial x}\bar{u}\bigr)\right) + 2\frac{\partial}{\partial x}\bigl(H_\delta\frac{\partial z_b}{\partial x}\bar{u}^2\bigr) + {\cal O}(\varepsilon^3\delta).
\end{eqnarray*}
In terms of the initial variables, the model (\ref{eq:sv_nonhyd_fin1})-(\ref{eq:sv_nonhyd_fin2}) corresponds to the one depicted in proposition~\ref{prop:SV++} with obvious expressions
for $H_m$, $H\bar{p}_{ng,nv}$ and $\left. p_{ng,nv}\right|_p$.

In order to obtain the energy equality for the model
(\ref{eq:sv++_1})-(\ref{eq:sv++_2}), we use the same
process and the same notations as in paragraph \ref{subsec:prop} but the
approximation order is now ${\cal O}(\varepsilon^3\delta)$ instead of
${\cal O}(\varepsilon^2\delta^2)$. Still using $\bar{u}= u + {\cal
O}(\varepsilon)= u_b + {\cal O}(\varepsilon)$, we have
\begin{eqnarray*}
{\tilde{\cal C}}_1 & = & \bar{u}\frac{\partial}{\partial x}\int_{z_b}^{\delta\eta} \Delta p + \bar{u}\frac{\partial z_b}{\partial x} \left. \Delta p\right|_b\nonumber\\
& = & \frac{\partial}{\partial x}\left(\int_{z_b}^{\delta\eta} u\Delta p \right) + [w \Delta p]_{z_b}^{\delta\eta} - \int_{z_b}^{\delta\eta} w \frac{\partial \Delta p }{\partial z}  + u_b\frac{\partial z_b}{\partial x} \left. \Delta p\right|_b + {\cal O}(\varepsilon^3\delta),\nonumber\\
\end{eqnarray*}
with $\Delta p = p - p_{nh}$, $p$ being given by (\ref{eq:pfin}). From
Eqs.~(\ref{eq:pfin}), (\ref{eq:ph}) and the boundary condition (\ref{eq:BC_ap1}), we get
$$\left.\Delta p\right|_s = {\cal O}(\varepsilon^3\delta),\quad \left.\Delta
p\right|_b = \varepsilon^2\delta^2\int_{z_b}^{\delta\eta}
\frac{\partial (uw)}{\partial x} + \varepsilon^2\delta^2 (w_s^2 - w_b^2)+ {\cal O}(\varepsilon^3\delta),$$
$$\frac{\partial \Delta p }{\partial z} = - \varepsilon^2\delta^2 \frac{\partial
(uw)}{\partial x} - 2\varepsilon^2\delta^2 w\frac{\partial
w}{\partial z} + {\cal O}(\varepsilon^3\delta),$$
leading to
\begin{eqnarray*}
{\tilde{\cal C}}_1 & = & \frac{\partial}{\partial x}\int_{z_b}^{\delta\eta} u\Delta p + \varepsilon^2\delta^2\int_{z_b}^{\delta\eta} w \frac{\partial uw}{\partial x} + \frac{2}{3}\varepsilon^2\delta^2 (w_s^3 - w_b^2)\nonumber\\
& & - \frac{\partial z_b}{\partial t}\left. \Delta p \right|_b + {\cal O}(\varepsilon^3\delta),\nonumber\\
& = & \frac{\partial}{\partial x}\left(\int_{z_b}^{\delta\eta} u\bigl( \Delta p + \frac{w^2}{2}\bigr)\right) - \frac{\partial z_b}{\partial t}\left. \Delta p \right|_b + {\cal O}(\varepsilon^3\delta).
\end{eqnarray*}
Returning to the initial variables, the preceding relation and the
expression of ${\cal R}$ obtained in paragraph~\ref{subsec:prop} allows us to write an energy equality for the model (\ref{eq:sv++_1})-(\ref{eq:sv++_2}) under the form
\begin{eqnarray}
\lefteqn{\frac{\partial \bar{E}}{\partial t} + \frac{\partial}{\partial x}
\left(\bar{u}\left( \bar{E} + H \bar{p}\right) - \nu\int_{z_b}^{\eta} \left(2H\bar{u}\frac{\partial u}{\partial x} + w\frac{\partial u}{\partial z}\right) + \kappa\left(\frac{H^2}{6}\frac{\partial \bar{u}}{\partial x}\bar{u} - \frac{H}{2}\frac{\partial z_b}{\partial x}\bar{u}^2\right)\right) \nonumber}\\
& & = - 2\nu\int_{z_b}^{\eta} \left(\left(\frac{\partial \bar{u}}{\partial x}\right)^2 + \left(\frac{\partial u}{\partial x}\right)^2\right)  - \frac{\kappa}{6}\left( H\frac{\partial \bar{u}}{\partial x} + \frac{\partial H}{\partial x}\bar{u}\right)^2\nonumber\\
& & - \frac{\kappa}{3}\left(\left(\frac{\partial z_b}{\partial x}- \frac{1}{4}\frac{\partial \eta}{\partial x}\right)^2 - \frac{1}{8}\left(\frac{\partial \eta}{\partial x}\right)^2\right)\bar{u}^2 -\frac{\kappa}{1+\frac{\kappa H}{3 \nu}}\left(1+\frac{3}{2}\left(\frac{\partial z_b}{\partial x}\right)^2\right) \bar{u}^2\nonumber\\
& & - H\frac{\partial p^a}{\partial t} + \left(\left. p_{nh}\right|_b + 2\nu\frac{\partial u_b}{\partial x}\right)\frac{\partial z_b}{\partial t} - \kappa\frac{H}{2}\frac{\partial^2 z_b}{\partial x\partial t}\bar{u},\nonumber
\end{eqnarray}
with
\begin{eqnarray*}
\bar{E} & = & \frac{H\overline{u^2}}{2} + \frac{H\overline{w^2}}{2} + \frac{gH(\eta+z_b)}{2},\nonumber\\
H \bar{p} & = & \int_{z_b}^{\eta} p \ dz,\quad H\overline{u^2} = H\left(1+\frac{2\kappa^2 H^2}{15\nu^2}\right)\bar{u}^2,\nonumber\\
H\overline{w^2} & = & \int_{z_b}^{\eta} w^2 =
H\left(\frac{\eta^2 + \eta z_b + z_b^2}{3}\left(\frac{\partial \bar{u}}{\partial
x}\right)^2 - (\eta + z_b) \frac{\partial \bar{u}}{\partial x}\frac{\partial
(z_b\bar{u})}{\partial x}\right.\\
& & \left.  + \left(\frac{\partial (z_b\bar{u})}{\partial
x}\right)^2 \right) + H\left(\frac{\partial z_b}{\partial t}\right)^2 + 2\frac{\partial z_b}{\partial t}\left(-\frac{\eta^2 - z_b^2}{2}\frac{\partial \bar{u}}{\partial x} + H \frac{\partial (z_b\bar{u})}{\partial x}\right) 
\end{eqnarray*}
Note that except for the friction terms, the previous expression is
analogous to the energy equality for the Navier-Stokes system
\cite{lions} but expressed with the vertically averaged variables. When the time derivatives of $p^a$ and $z_b$ are dropped, the right hand
side of the preceding energy equality is negative when $\frac{\partial \eta}{\partial x}$ is enough small.

\section{Conclusion}

In this paper we have derived two extensions of the Saint-Venant system when
the hydrostatic assumption is relaxed. The obtained models, especially
in section \ref{sec:SV_nhyd1}, are similar to Boussinesq type models but
derived in a more rigourous context and satisfying an energy equality.

On one hand the averaged models of shallow water type presented in
this paper reduce the complexity of the discretization of the Navier-Stokes
equations since they are written over a fixed domain. But on the other hand their mathematical formulation is
more complex since high order derivatives~-- especially in space~-- appear.

The preliminary numerical simulations and comparison with experimental
measurements performed with the proposed models are promising. They are
not presented in this paper and will be described in a forthcoming publication.

\bigskip
{\bf Acknowledgements.} The authors want to thank Emmanuel Audusse, Fran\c{c}ois Bouchut and Beno\^{\i}t Perthame for helpful discussions that have allowed to greatly improve the paper.

\bibliographystyle{amsplain}
\bibliography{boussinesq}

\providecommand{\bysame}{\leavevmode\hbox to3em{\hrulefill}\thinspace}
\providecommand{\MR}{\relax\ifhmode\unskip\space\fi MR }
\providecommand{\MRhref}[2]{%
  \href{http://www.ams.org/mathscinet-getitem?mr=#1}{#2}
}
\providecommand{\href}[2]{#2}
\begin{thebibliography}{10}

\bibitem{arnold}
D.N. Arnold, J.~Douglas, and V.~Thom\'ee, \emph{{Superconvergence of a Finite
  Element Approximation to the Solution of a Sobolev Equation in a Single Space
  Variable}}, Mathematics of Computation \textbf{36} (1981), no.~153, 53--64.

\bibitem{audusse}
E.~Audusse, \emph{{A multilayer Saint-Venant System~: Derivation and Numerical
  Validation}}, Discrete and Continuous Dynamical Systems, Ser. B \textbf{5}
  (2005), no.~2, 189--214.

\bibitem{saint-venant}
A.J.C. Barr\'e~de Saint-Venant, \emph{{Th\'eorie du mouvement non permanent des
  eaux avec applications aux crues des rivi\`eres et \`a l'introduction des
  mar\'ees dans leur lit}}, C. R. Acad. Sci. Paris \textbf{73} (1871),
  147--154.

\bibitem{bbm}
J.L. Bona, T.B. Benjamin, and J.J. Mahony, \emph{{Model equations for long
  waves in nonlinear dispersive systems}}, Philos. Trans. Royal Soc. London
  Series A \textbf{272} (1972), 47--78.

\bibitem{saut1}
J.L. Bona, M.~Chen, and J.C. Saut, \emph{{Boussinesq equations and other
  systems for small-amplitude long waves in nonlinear dispersive media: Part
  {I}. Derivation and linear theory}}, J. Nonlinear Sci. \textbf{12} (2002),
  283--318.

\bibitem{saut2}
\bysame, \emph{{Boussinesq equations and other systems for small-amplitude long
  waves in nonlinear dispersive media: Part {II}. Nonlinear theory}},
  Nonlinearity \textbf{17} (2004), 925--952.

\bibitem{bouchut}
F.~Bouchut and M.~Westdickenberg, \emph{{Gravity driven shallow water models
  for arbitrary topography}}, Comm. in Math. Sci. \textbf{2} (2004), 359--389.

\bibitem{boussinesq1}
J.V. Boussinesq, \emph{{Th\'eorie de l'intumescence liquide appel\'ee onde
  solitaire ou de translation se propageant dans un canal rectangulaire}}, C.
  R. Acad. Sci. Paris \textbf{72} (1871), 755--759.

\bibitem{boussinesq2}
\bysame, \emph{{Th\'eorie g\'en\'erale des mouvements qui sont propag\'es dans
  un canal rectangulaire horizontal}}, C. R. Acad. Sci. Paris \textbf{73}
  (1871), 256--260.

\bibitem{boussinesq3}
\bysame, \emph{{Th\'eorie des ondes et des remous qui se propagent le long d'un
  canal rectangulaire horizontal, en communiquant au liquide contenu dans ce
  canal des vitesses sensiblement pareilles de la surface au fond}}, J. Math.
  Pures Appl. \textbf{17} (1872), 55--108.

\bibitem{cienfuegos1}
R.~Cienfuegos, E.~Barth\'elemy, and P.~Bonneton, \emph{{A fourth-order compact
  finite volume scheme for fully nonlinear and weakly dispersive
  Boussinesq-type equations. Part I: Model development and analysis}}, Int. J.
  Numer. Meth. Fluids \textbf{51} (2006), no.~11, 1217--1253.

\bibitem{cienfuegos2}
\bysame, \emph{{A fourth-order compact finite volume scheme for fully nonlinear
  and weakly dispersive Boussinesq-type equations. Part II: Boundary conditions
  and validation}}, Int. J. Numer. Meth. Fluids \textbf{53} (2006), no.~9,
  1423--1455.

\bibitem{saleri}
S.~Ferrari and F.~Saleri, \emph{{A new two-dimensional Shallow Water model
  including pressure effects and slow varying bottom topography}}, M2AN
  \textbf{38} (2004), no.~2, 211--234.

\bibitem{gerbeau}
J.-F. Gerbeau and B.~Perthame, \emph{{Derivation of Viscous Saint-Venant System
  for Laminar Shallow Water; Numerical Validation}}, Discrete and Continuous
  Dynamical Systems, Ser. B \textbf{1} (2001), no.~1, 89--102.

\bibitem{levermore}
C.D. Levermore and M.~Sammartino, \emph{A shallow water model with eddy
  viscosity for basins with varying bottom topography}, Nonlinearity
  \textbf{14} (2001), no.~6, 1493--1515.

\bibitem{lions}
P.L. Lions, \emph{{Mathematical Topics in Fluid Mechanics. Vol.~1:
  Incompressible models.}}, Oxford University Press, 1996.

\bibitem{marche}
F.~Marche, \emph{{Derivation of a new two-dimensional viscous shallow water
  model with varying topography, bottom friction and capillary effects}},
  European Journal of Mechanic /B \textbf{26} (2007), 49--63.

\bibitem{valentin}
B.~Mohammadi, O.~Pironneau, and F.~Valentin, \emph{{Rough boundaries and wall
  laws}}, Int. J. Numer. Meth. Fluids \textbf{27} (1998), no.~1-4, 169--177.

\bibitem{nwogu}
O.~Nwogu, \emph{{Alternative form of Boussinesq equations for nearshore wave
  propagation}}, Journal of Waterway, Port, Coastal and Ocean Engineering, ASCE
  \textbf{119} (1993), no.~6, 618--638.

\bibitem{peregrine}
D.H. Peregrine, \emph{{Long waves on a beach}}, J. Fluid Mech. \textbf{27}
  (1967), 815--827.

\bibitem{perotto}
S.~Perotto and F.~Saleri, \emph{{Adaptive finite element methods for Boussinesq
  equations}}, Numer. Methods Partial Differential Equations \textbf{16}
  (2000), no.~2, 214--236.

\bibitem{zech}
S.~Soares~Frazao and Y.~Zech, \emph{{Undular bores and secondary waves -
  Experiments and hybrid finite-volume modelling}}, Journal of Hydraulic
  Research \textbf{40} (2002), no.~1, 33--43.

\bibitem{ursell}
F.~Ursell, \emph{{The long wave paradox in the theory of gavity waves}}, Proc.
  Cambridge Phil. Soc. \textbf{49} (1953), 685--694.

\bibitem{walkley}
M.A. Walkley, \emph{{A numerical Method for Extended Boussinesq Shallow-Water
  Wave Equations}}, Ph.D. thesis, University of Leeds, 1999.

\end{thebibliography}


\end{document}